\documentclass[a4wide, oneside]{article}  

\usepackage{myPackages}
\usepackage{myAssignments3}
\usepackage{myOperators}

\newcommand{\PD}{\operatorname{PD}}  

\renewcommand{\S}{\Sigma}

\usepackage{tikz}
\usetikzlibrary{arrows,decorations.pathmorphing,backgrounds,positioning,fit,petri}
\usetikzlibrary{calc,arrows.meta,positioning}
\usepackage{ragged2e}

\usepackage{fancyhdr,floatpag}
\floatpagestyle{fancy}

\usepackage{caption}
\usepackage{subcaption}

\titleformat{\section}
{\normalfont\normalsize\bfseries\center}{\thesection}{1em}{}

\titleformat{\subsection}
{\normalfont\normalsize\bfseries}{\thesection}{1em}{}

\begin{document}
\renewcommand{\abstractname}{\vspace{-\baselineskip}}

\title{Topological constraints for Stein fillings of tight structures on lens spaces}

\author{Edoardo Fossati}

\date{}

\newcommand{\Addresses}{{
  \bigskip

  E.~Fossati, \textsc{Scuola Normale Superiore, Piazza dei Cavalieri 7, Pisa.}\par\nopagebreak
  \textit{E-mail address:} \texttt{edoardo.fossati@sns.it}

}}

\maketitle

\begin{abstract}
In this article we give a sharp upper bound on the possible values of the Euler characteristic for a minimal symplectic filling of a tight contact structure on a lens space. This estimate is obtained by looking at the topology of the spaces involved, extending this way what we already knew from the universally tight case to the virtually overtwisted one. As a lower bound, we prove that virtually overtwisted structures on a certain family lens spaces never bound Stein rational homology balls.

\noindent Then we turn our attention to covering maps: since an overtwisted disk lifts to an overtwisted disk, all the coverings of a universally tight structure are themselves tight. The situation is less clear when we consider virtually overtwisted structures. By starting with such a structure on a lens space, we know that this lifts to an overtwisted structure on $S^3$, but what happens to all the other intermediate coverings? We give necessary conditions for these lifts to still be tight, and deduce some information about the fundamental groups of the possible Stein fillings of certain virtually overtwisted structures.
\end{abstract}

\section{Introduction}

Classifying symplectic fillings (up to homeomorphism, diffeomorphism or symplectic deformation equivalence) of a given contact 3-manifold can be a very hard task, even though some progress has been made in the last years. Lens spaces surely represent a class of 3-manifolds for which many results are known. McDuff showed in \cite{mcduff} that $L(p,1)$, endowed with the standard tight contact structure, has a unique Stein filling when $p\neq 4$, and two different Stein fillings when $p=4$. Later, Lisca \cite{lisca} extended McDuff's results and provided a procedure to produce all of the Stein fillings of $(L(p,q),\xi_{st})$, up to diffeomorphism. Partial results about fillings are available when one considers non-standard tight contact structures on lens spaces, i.e. those that pull back to an overtwisted structure on the universal cover $S^3$ and that are therefore called \emph{virtually overtwisted} (see \cite{plamenVHM}, \cite{kaloti}, \cite{fossati}). A more modest approach is trying to give some constraints on the topological invariants of the Stein fillings, even if a complete classification is missing. 

If we restrict to planar contact structures, then studying Stein fillings is enough if we want to understand weak symplectic fillings, since these are symplectically deformation equivalent to blow ups of Stein fillings, see \cite{niederkruger}.
Some topological constraints for Stein fillings of planar contact structures have already been found (see for example \cite{etnyre}, \cite{c1planar}, \cite{wendl}, \cite{wand}).

In this article we focus our attention on lens spaces. Let $p>q>1$ be coprime integers and compute the continued fraction expansion
\[\frac{p}{q}=[a_1,a_2,\ldots, a_l]=a_1-\frac{1}{a_2-\frac{1}{ \mbox{\tiny $\ddots$}\, -\frac{1}{a_l}}},\]
with $a_i\geq 2\;\forall i=1,\ldots,l$. We will often refer to $l$ as the length of this expansion.
To this expansion we can associate a negative linear graph $\Lambda(p,q)$ and a corresponding negative definite 4-manifold $X_{\Lambda(p,q)}$ realized as a plumbing. We give a sharp upper bound on the possible values of the Euler characteristic for a minimal symplectic filling of a tight contact structure on a lens space:

\begin{thm} \label{maxchi}
Let $\xi$ be any tight contact structure on $L(p,q)$. Let $W$ be a minimal symplectic filling of $L(p,q)$ and let $l=\length(p/q)$. Then
\[\chi(W)\leq 1+l.\]
\end{thm}

This is just a topological constraint, not involving contact structures, and therefore it is valid both for the universally tight and virtually overtwisted cases. Moreover, the upper bound is always realized by a minimal symplectic filling ($X_{\Lambda(p,q)}$ itself supports a Stein structure inducing the prescribed contact structure on its boundary) whose intersection form and fundamental group are uniquely determined:

\begin{thm}\label{uniqueintform}
Let $\xi$ be any tight contact structure on $L(p,q)$ and let $l=\length(p/q)$. Let $X$ be a minimal symplectic filling of $(L(p,q),\xi)$ with $b_2(X)=l$. Then the intersection form $Q_X$ is isomorphic to the intersection form of $X_{\Lambda(p,q)}$. Moreover, $X$ is simply connected.
\end{thm}

We also prove the following corollary, regarding the uniqueness (in certain cases) of the filling with maximal Euler characteristic:

\begin{cor}\label{homeo}
Let $\xi$ be a tight contact structure on $L(p,q)$ and let $l=\length(p/q)$. Let $X$ be a minimal symplectic filling of $(L(p,q),\xi)$ with $b_2(X)=l$. Assume that $p\in\{2,4,s^n,2s^n\}$, for any odd prime $s$ and any $n$. Then $X$ is homeomorphic to $X_{\Lambda(p,q)}$.
\end{cor}

On the other hand, the first and third Betti numbers of a Stein filling $W$ of a lens space $L(p,q)$ are always zero (\cite[page 216]{ozbagci}), hence we have an obvious lower bound on the value $\chi(W)$, which is $\chi(W)=b_0(W)=1$: this is realized precisely when $(L(p,q),\xi)$ bounds a Stein rational homology ball. In \cite{liscaribbon} and \cite[page 247]{ana} it is proved that, in order to guarantee the existence of rational balls with boundary a lens space, the numbers $p$ and $q$ must fall into one of three families with specific numerical conditions.

Among those, we restrict to the case when $p$ and $q$ are of the form $p=m^2$ and $q=mk-1$, for some $m>k>0$ with $(m,k)=1$. It is known (\cite[Corollay 1.2c]{lisca}) that $L(m^2,mk-1)$ endowed with a universally tight contact structure bounds a Stein rational ball and we use this fact to prove that in the virtually overtwisted case this never happens, concluding that:

\begin{thm} \label{minchi}
Let $W$ be a symplectic filling of $(L(p,q),\xi_{vo})$, with $p=m^2$ and $q=mk-1$, for some $m>k>0$ and $(m,k)=1$. Then $\chi(W)\geq 2$.
\end{thm}

Theorem \ref{minchi} can be generalized to the other families of lens spaces which are known to bound a smooth rational homology ball: these balls do not support any symplectic structure, i.e. none of the virtually overtwisted contact structures can be filled by a Stein rational ball, see \cite[Proposition A.1]{gollastar}.

In Section \ref{coveringsection} we give a series of examples using the description of tight structures given by Honda in \cite{honda} to study explicit cases of covering maps between lens spaces endowed with contact structures. One of the problems we faced when studying such covering is the mysterious behavior of numbers: for example, there is no understanding on how the lengths of $p/q$ and $p'/q$ are related, if $p'$ is a divisor of $p$. This makes the problem hard even to organize, since we could not glimpse any clear scheme or pattern for stating reasonable guesses. Theorem \ref{thmcovering} is the only stance of a general result which does not depend on specific examples.

\begin{thm}\label{thmcovering}
Let $p,\,q$ and $d$ be such that $q<p<dq$. Then every virtually overtwisted contact structure on $L(p,q)$ lifts along a degree $d$ covering to a structure which is overtwisted.
\end{thm}

The last part is dedicated to the study of the fundamental group of Stein fillings of virtually overtwisted structures on lens spaces, combining the results above about Euler characteristic with what we developed on the behavior of coverings. Recall that the fundamental group of any Stein filling is a quotient of the fundamental group of its boundary (see \cite{os}). 
As a consequence of Theorem \ref{maxchi}, we will prove Theorem \ref{pi1b2} and provide some specific examples and applications.

\begin{thm}\label{pi1b2} 
Let $X$ be a Stein filling of $(L(p,q),\xi)$ with $\pi_1(X)=\Z/d\Z$, for $p=dp'$. Then
\[\chi(X)\leq \frac{1+l'}{d},\]
where $l'=\length(p'/q')$, with $q'\equiv q \pmod{p'}$.
\end{thm}

\paragraph{Acknowledgments} The author wishes to thank Paolo Lisca for all the useful discussions and comments on the contents of this article, which will be part of his PhD thesis. Thanks are also due to Ko Honda for helpful email correspondence. 
A few remarks on the first version of the paper came from Marco Golla, who found a gap in the proof of Proposition \ref{pi=1} (which has been filled now) and indicated further unclear points. Moreover, he provided examples of non simply-connected Stein fillings of some virtually overtiwsted structures on certain lens spaces, answering a question that was originally formulated in the first version of this paper and whose solution is now presented at the end of the work. Lastly, the author thanks the referee for the interesting comments and remarks.

\section{Upper bound for the Euler characteristic}

The goal of this section is to prove Theorems \ref{maxchi} and \ref{uniqueintform}. First, recall that a vertex $v$ of a weighted graph is a \emph{bad vertex} if 
\[w(v)+d(v)>0,\]
where $w(v)$ and $d(v)$ are respectively the weight and the degree (i.e. the number of edges containing $v$) of the vertex. Then:

\begin{thm}\label{maxtree}
Let $\Gamma$ be a negative definite plumbing tree with $k$ vertices, none of which is a bad vertex. Call $\overline{Y}$ the plumbed 3-manifold associated to $\Gamma$ and assume that $\overline{Y}$ is a rational homology sphere. Denote by $Y$ the manifold with the opposite orientation. Let $X$ be a negative definite smooth 4-manifold with no $(-1)$-class in $H_2(X;\Z)$ such that $\partial X=Y$. Then
\[b_2(X)\leq 1+\sum_{i=1}^k (|w(v_i)|-2).\]
\end{thm}

\begin{prf} 
Let $P=P_{\Gamma}$ the plumbed 4-manifold associated to $\Gamma$, whose oriented boundary is $\overline{Y}$, and whose intersection form is $Q_P=Q_{\Gamma}$. Form the closed manifold $W=X\cup_{\partial}P$ by gluing the two manifolds along the boundary, see Figure \ref{glue1}. 
\begin{figure}[ht!]
	\center
	\includegraphics[scale=0.5]{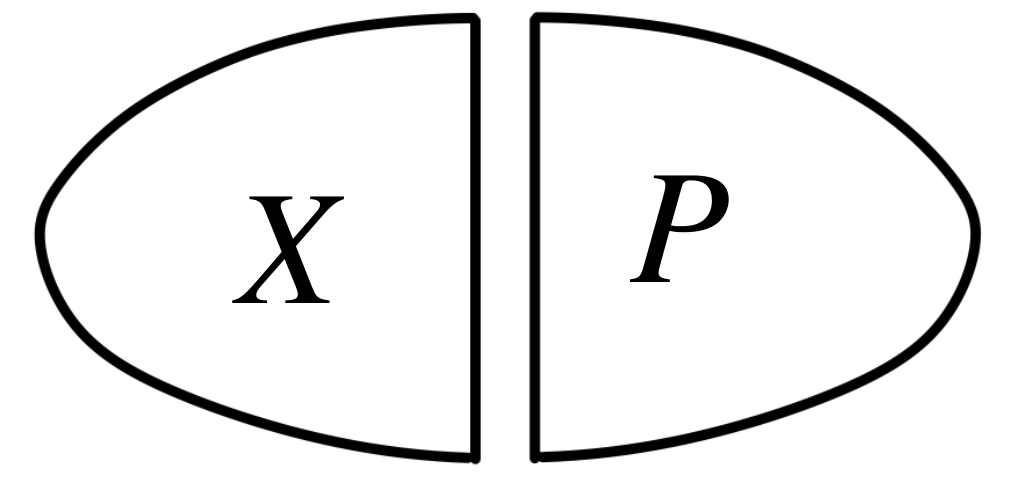}
		\caption{The closed manifold $W=X\cup_{\partial}P$.}
	\label{glue1}
\end{figure}

\noindent We get a closed smooth 4-manifold $W$ whose intersection form $Q_W$ is negative definite, and hence, by Donaldson's theorem \cite{donaldson}, isomorphic to $\langle -1\rangle^r$, for some $r$. Since $Y$ is a rational homology sphere, we have that
\[\rank(Q_X)+\rank(Q_{\Gamma})=b_2(X)+b_2(P)=b_2(X)+k=r.\]
We want to show that 
\[Q_X\simeq (Q_{\Gamma})^{\perp}.\]
A priori, $Q_X$ is a sub-lattice of finite index $n$, for some $n>0$:
\[Q_X\leq (Q_{\Gamma})^{\perp},\]
but the following argument (between the star symbols) shows that they coincide. 

\noindent $\color{blue}{\star}$ Look at the exact sequence of the pair $(W,X)$:
\[ H_3(W,X)\to H_2(X) \to H_2(W) \to H_2(W,X) \]
and notice that, by excision and Poincaré-Lefschetz duality
\[H_3(W,X)\simeq H_3(P,\partial P)\simeq H^1(P)=0\] 
and similarly
\[H_2(W,X)\simeq H_2(P,\partial P)\simeq H^2(P).\]
The latter is free, because its torsion comes from $H_1(P)$, which is $0$. Therefore, the inclusion $H_2(X)\hookrightarrow H_2(W)$ has a free quotient, being this a subgroup of the free group $H_2(W,X)$. Hence, if we take a class $\alpha\in H_2(W)$ with the property that $n\alpha$ is inside $H_2(X)$, we automatically get $\alpha\in H_2(X)$. In particular, if $\alpha\in (Q_{\Gamma})^{\perp}$, then $n\alpha\in (H_2(X),Q_X)$, with $n$ equal to the index $Q_X\leq (Q_{\Gamma})^{\perp}$, and hence $\alpha\in (H_2(X),Q_X)$. So $Q_X\simeq (Q_{\Gamma})^{\perp}$, as wanted. $\color{blue}{\star}$

The isomorphism 
\[Q_X\simeq (Q_{\Gamma})^{\perp}\]
implies that we have an embedding $\Gamma\hookrightarrow \langle -1\rangle^r$ with the property that there is no $(-1)$-class in the orthogonal, otherwise this would come from $X$, which, by assumption, does not have any. We call such an embedding \emph{irreducible}. 
Notice that, up to isomorphism, there is a unique maximal irreducible embedding
\[\Gamma\hookrightarrow \langle -1\rangle^t,\]
where by \emph{maximal} we mean that the dimension $t$ (a priori infinite) of the ambient lattice cannot be bigger. First of all notice that at least one irreducible embedding exists: we will explicitly describe the construction of one of them, which turns out to be the maximal one. Since the sum of the weights of the graph is finite, $t$ is itself finite.
To embed the graph $\Gamma$ in such a way that there is no $(-1)$-class in the orthogonal complement implies that all the elements in the canonical basis $\{e_1,\ldots,e_t\}$ of $\langle -1\rangle^t$ appear in the image of some vertex of the graph. Therefore, to obtain the maximal such $t$, we have to impose only the requirements that:
\begin{itemize}
\item[1)] the $i^{th}$ vertex is sent to a combination of $|w_i|$-many distinct basis elements and
\item[2)] any two adjacent vertices of $\Gamma$ share, via the embedding, exactly one element $e_j$.
\end{itemize}
If one of these conditions is not satisfied, then we end up with (at least) one line $\langle e_j \rangle$ which is not hit by the image of $\Gamma$ and that will produce an element in the orthogonal with square $-1$. So the image of the first vertex with weight $w_1$ must be a sum of $|w_1|$-many distinct elements $e_i$. The second vertex is sent to a combination of $|w_2|$ elements, among which exactly one has already appeared in the image of the first vertex, and so on. Hence we find
\[t=1+\sum_{i=1}^k (|w(v_i)|-1).\]
Therefore, since the dimension of the maximal irreducible embedding of $\Gamma$ is as above, we have $r=k+b_2(X)\leq t$. We conclude:
\[b_2(X)\leq  t-k=1+\sum_{i=1}^k (|w(v_i)|-1)-k=1+\sum_{i=1}^k (|w(v_i)|-2).\]
\end{prf}

\begin{cor}\label{uniquenessgeneral}
In the setting above, the intersection form of the the manifolds with boundary $Y$ and maximal $b_2$ is uniquely determined, up to isomoprhism.
\end{cor}

\begin{prf} Assume that $X_1$, $X_2$ are negative definite with no $(-1)$-class, $\partial X_1=\partial X_2=Y$ and with $b_2(X_i)$ maximal. Then, by uniqueness of the maximal irreducible embedding $\Gamma\hookrightarrow \langle -1\rangle^t$, we have that 
\[Q_{X_1}\simeq Q_{X_2}\simeq \Gamma^{\perp}\subseteq \langle -1\rangle^t.\]
\end{prf}

Now we specialize to the case of lens spaces. Start with $L(p,q)$ and take the expansion $p/q=[a_1,a_2,\ldots, a_l]$, where all the $a_i$'s are $\geq 2$. Call $\Lambda=\Lambda(p,q)$ the associated negative definite lattice with $l$ vertices (where $l=\length(p/q)$):
\begin{center}
$\Lambda=$
\begin{tikzpicture}
        \node[shape=circle,fill=black,inner sep=1.5pt,label=$-a_1$] (1)                  {};
        \node[shape=circle,fill=black,inner sep=1.5pt,label=$-a_2$] (2) [right=of 1] {}
        edge [-]               (1);
        \node[shape=circle,fill=black,inner sep=1.5pt,label=$-a_3$] (3) [right=of 2] {}
        edge [-]               (2);
        \node[shape=circle,fill=black,inner sep=1.5pt,label=$-a_{l-1}$] (4) [right=of 3] {} ;
                 \node at ($(3)!.5!(4)$) {\ldots};

        \node[shape=circle,fill=black,inner sep=1.5pt,label=$-a_l$] (5) [right=of 4] {}
        edge [-]               (4);
\end{tikzpicture}
\end{center}
We apply Riemenschneider's dots method \cite{riemen} to build a negative definite 4-manifold with boundary of $L(p,p-q)$ whose intersection lattice will be called $\Lambda^{\nu}$. This is obtained by reading column-wise the entries of Table \ref{table}.
\begin{center}
\begin{table}[h!]
        \centering
        \begin{tabular}{ cccccccccccc}
            \multicolumn{4}{l}{$\overbrace{\rule{2.2cm}{0pt}}^{a_1-1}$}  &  & & & & & & & \\
           $\bullet$  & $\bullet$ &  $\cdots$ &  $\bullet$  &   &  &   &   &   &  &   &    \\
             \hline 
           & & &$\bullet$  & $\bullet$ &  $\cdots$ &  $\bullet$  &   &  &   &   &   \\
             \hline 
           & & &  & &   &    & $\ddots$   &  &   &   &   \\
             \hline 
           & & &  & &   &    &   &   $\bullet$  & $\bullet$ &  $\cdots$ &  $\bullet$\\
           & & &  \multicolumn{4}{l}{$\underbrace{\rule{2.2cm}{0pt}}_{a_2-1}$} & &  \multicolumn{4}{l}{$\underbrace{\rule{2.2cm}{0pt}}_{a_l-1}$} 
 \end{tabular}
  \caption{Riemenschneider's dots method.}
  \label{table}
  \end{table}
\end{center}

\noindent If we call $l^{\nu}$ the number of columns and set 
\[c_j=1+\#\{\mbox{dots in the $j^{\mbox{\tiny{th}}}$ column}\},\]
then we obtain the continued fraction expansion of $p/(p-q)$ as
\[\frac{p}{p-q}=[c_1,c_2,\ldots,c_{l^{\nu}}].\]
Before proving Theorem \ref{maxchi} that, we need a lemma.

\begin{lem}\label{lem1}
\[\length(p/q)+\length(p/(p-q))=1+\sum_{i=1}^l(a_i-1).\]
\end{lem}

\begin{prf}
We know that $\length(p/q)=l$, so we compute $\length(p/(p-q))$. This is just the number $l^{\nu}$ of columns:
\begin{align*}
l^{\nu} = & (a_1-1)+(a_2-2)+\ldots+(a_l-2)\\
=\, &\sum_{i=1}^l(a_i-2)+1\\
=\, &\sum_{i=1}^l(a_i-1)-l+1.
\end{align*}
Therefore $\length(p/q)+\length(p/(p-q))=1+\sum_{i=1}^l(a_i-1)$. 
\end{prf}

\noindent By switching the roles of $p/q$ and $p/(p-q)$, it is clear from Lemma \ref{lem1} that
\begin{align*}
\length(p/q)+\length(p/(p-q))=& \rank(\Lambda)+\rank(\Lambda^{\nu})\\
=\,& l+l^{\nu}\\
=\,& 1+\sum_{i=1}^l(a_i-1)\\
=\, & 1+\sum_{i=1}^{l^{\nu}}(c_i-1).
\end{align*}

\begin{nrem}\label{b1zero}
In the book \cite[Section 12.3]{ozbagci} the authors made the following observation, which we will often use in this work. If $X$ is a Stein filling of $Y$, then the morphism $\pi_1(Y)\to \pi_1(X)$, induced by the inclusion, is surjective since $X$ can be built on $Y \times [0, 1]$ by attaching 2-, 3- and 4-handles only. In particular, $b_1(X)\leq b_1(Y)$ and if $Y$ is a lens space, then $b_1(X)=0$.
\end{nrem}

Theorem \ref{maxchi} follows now from Theorem \ref{maxtree}:

\begin{prf}[of Theorem \ref{maxchi}]
Let $Y=L(p,q)$, so that $\overline{Y}$ is the 3-manifold associated to $\Lambda^{\nu}$, with $\Lambda^{\nu}$ playing the role of $\Gamma$. The setting for lens spaces is coherent with the hypotheses of Theorem \ref{maxtree}:
\begin{itemize}
\item lens spaces arise as plumbings on trees with no bad vertices;
\item lens spaces are rational homology spheres;
\item contact structures on lens spaces are planar (\cite[Theorem 3.3]{schonenberger}), and therefore $b_2(X)=b_2^-(X)$ for any minimal filling $X$ (\cite{etnyre});
\item minimal fillings of planar contact structures have no $(-1)$-class, as proved in \cite[Corolloary 1.8]{ghigginigolla};
\end{itemize}

Therefore, since any minimal filling $X$ of $(Y,\xi)$ has $b_1=0$, as Remark \ref{b1zero} explained, we have:
\begin{align*}
\chi(X)= & 1+b_2(X)\\
\leq\, & 1+\left(1+\sum_{i=1}^{l^{\nu}} (c_i-2)\right)\\
=\, & 1+\left(1+\sum_{i=1}^{l^{\nu}} (c_i-1)-l^{\nu}\right)\\
=\, & 1+(l+l^{\nu}-l^{\nu})\\
=\, & 1+l.
\end{align*}
\end{prf}

\begin{prf}[of Theorem \ref{uniqueintform}]
The fact that the intersection form is uniquely determined is just a special case of Corollary \ref{uniquenessgeneral}. For the fundamental group, let $X$ be a filling with $b_2(X)=l=\length(p/q)$. We know that
\[Q_X\simeq Q_{X_{\Lambda(p,q)}},\]
and we look at the long exact sequence of the pair $(X,\partial X)$, with $\partial X=L(p,q)$:
\[\xymatrix{ H_2(L(p,q))\ar[r] & H_2(X)\ar[r]^-{Q_X} &H_2(X,L(p,q))\ar[r]& H_1(L(p,q))\ar[r]& H_1(X)\ar[r] &H_1(X,L(p,q))}\]
Notice that:
\begin{itemize}
\item[1)] $H_2(L(p,q))\simeq H^1(L(p,q))=0$;
\item[2)] $H_2(X))\simeq \Z^l$;
\item[3)] $H_2(X,L(p,q)\simeq H^2(X)\simeq \Z^l\oplus H_1(X)$;
\item[4)] $H_1(L(p,q))\simeq \Z/p\Z$;
\item[5)] $H_1(X,L(p,q))\simeq H^3(X)=0$;
\item[6)] $\det(Q_X)=p$.
\end{itemize}
Therefore, by substituting everything, it follows that $H_1(X)=0.$ But since, by Remark \ref{b1zero}, $\pi_1(X)$ is abelian, we have that $\pi_1(X)=0$, as wanted.
\end{prf}

\noindent We can now give a proof of Corollary \ref{homeo}.

\begin{prf}[of Corollary \ref{homeo}] To prove this corollary we need three facts:

\begin{itemize}
\item[1)] $X$ and $X_{\Lambda(p,q)}$ are both simply connected by Theorem \ref{uniqueintform};
\item[2)] $X$ and $X_{\Lambda(p,q)}$ have isomorphic intersection forms by Theorem \ref{uniqueintform};
\item[3)] the fundamental group of their boundary is $\pi_1(Y)\simeq \Z/p\Z$, with $p\in\{2,4,s^n,2s^n\}$, for any odd prime $s$ and any $n$.
\end{itemize}

\noindent Then \cite[Proposition 0.6]{boyer} applies and tells that $X$ and $X_{\Lambda(p,q)}$ are homeomorphic.
\end{prf}

\section{Lower bound for the Euler characteristic} \label{lowerboundsection}

The goal of this section is to prove Theorem \ref{minchi}, i.e. that, among the virtually overtwisted structures on the lens spaces of the form $L(m^2,mk-1)$ with $(m,k)=1$, none of these can be filled by a Stein rational homology ball. 

The first thing to notice is that, thanks to Honda's classification result \cite{honda}, each tight contact structure on a lens space has a Legendrian surgery presentation which comes from placing the corresponding chain of unknots into Legendrian position with respect to the standard contact structure of $S^3$. So, by looking at the rotation numbers of the various components of the link, we can describe all the tight contact structures that a lens space supports, up to isotopy.

Let $(Y,\xi)$ be a contact 3-manifold with $c_1(\xi)$ a torsion class. Then, in \cite{gompfsymp}, Gompf defined the invariant
\[d_3(Y,\xi)=\frac{1}{4}(c_1(X,J)^2-3\sigma(X)-2\chi(X))\in \Q,\]
where $(X,J)$ is any almost complex 4-manifold with boundary $\partial X=Y$ such that $\xi$ is homotopic to $TY\cap JTY$ (compare with Lemma 6.2.6 of \cite{ozbagci}).

\begin{lem}
If $(Y,\xi)$ bounds a Stein rational homology ball, then $d_3(Y,\xi)=-\frac{1}{2}$.
\end{lem}

\begin{prf} The quantity $d_3=\frac{1}{4}(c_1^2-3\sigma-2\chi)$ does not depend on the chosen filling, and if $(Y,\xi)=\partial (X,J)$ with $H_2(X;\Q)=H_1(X;\Q)=0$, then
\[d_3(Y,\xi)=\frac{1}{4}(c_1(X,J)^2-3\sigma(X)-2\chi(X))=\frac{1}{4}(0-0-2)=-\frac{1}{2}.\]
\end{prf}

\noindent In the case of lens spaces, the computation of the $d_3$ invariant is as follows:
\[d_3=\frac{1}{4}(c_1^2-3	\sigma-2(1-\sigma))=\frac{1}{4}(c_1^2-\sigma-2),\]
because all the Stein fillings have $b_1=b_3=0$ (see Remark \ref{b1zero}) and $b_2=b_2^-$ (\cite{etnyre}).
This means that, if $(L(p,q),\xi)$ bounds a Stein rational ball, then for \emph{any} other filling $(X,J)$ we have:
\[-\frac{1}{2}=\frac{1}{4}(c_1(J)^2-\sigma(X)-2)\]
and hence
\begin{equation}\label{c2=sigma}
c_1(J)^2=\sigma(X).
\end{equation}

We want to compute $c_1(J)^2$ for the filling $(X,J)$ of $(L(p,q),\xi)$ realized as plumbing described by the linear graph of the expansion $p/q$.
To do this, we need to specify the vector $r$ of rotation numbers for the components of the linear plumbing. If
\[\frac{p}{q}=[v_1,v_2,\ldots v_n]=v_1-\frac{1}{v_2-\frac{1}{\ddots -\frac{1}{v_n}}},\]
with all $v_i\geq 2$, then the quantity $c_1(J)^2$ is given by
\begin{equation}\label{c2=rot}
r^T(Q)^{-1}r,
\end{equation}
where $Q$ is the matrix
\[Q=
\begin{bmatrix}
    -v_1       & 1  \\
 1       & \ddots & \ddots \\
 & \ddots & \ddots & \ddots \\
  & & \ddots &\ddots &1\\
  & & & 1 & -v_n\\
\end{bmatrix}
,\]
which represents the intersection form of $X$ in the basis corresponding to the linear graph, where each vertex is a generator.
According to \cite{honda}, there are two universally tight contact structure on $L(p,q)$ up to isotopy (and just one on $L(p,p-1)$). Honda also characterizes the rotation number of each component of the link given by the chain of Legendrian unknots, whose associated Legendrian surgered manifold is $(L(p,q),	\xi_{st})$.

Let $y=(-v_1+2,-v_2+2,\ldots,-v_n+2)$ be the vector of these rotation numbers, i.e. the vector corresponding to one of the two universally tight (standard) structures on $L(p,q)$, the other one being $-y$. By construction, the rotation vectors representing the virtually overtwisted structures have components $x_i$ satisfying 
\[|x_i|\leq |y_i|,\]
with at least one index $\overline{\imath}$ for which $|x_{\overline{\imath}}|<|y_{\overline{\imath}}|$.

\noindent Consider the function $f:\R^n\to \R$ given by $z\mapsto \|z\|_{Q^{-1}}=z^T(Q)^{-1}z$ and notice that, by Equalities \eqref{c2=rot} and \eqref{c2=sigma},
\[f(y)=\sigma(P)=-n.\]
Theorem \ref{minchi} follows directly from Proposition \ref{concave} below, but first we need:

\begin{lem} \label{matrix}
All the entries of the matrix $Q^{-1}$ are strictly negative (in short: $Q^{-1}\ll 0$).
\end{lem}

\begin{prf}
The condition $Q^{-1}\ll 0$ is true if we show that $Q^{-1}x\ll 0$ holds whenever $x$ is a non-zero vector with non-negative components, i.e. $0\neq x \underline{\gg} 0$ (this is just a consequence of the fact that the columns of $Q^{-1}$ are the images of the vectors of the canonical basis).

So we need to check that: $0\neq x\;\underline{\gg}\; 0$ implies $Q^{-1}x\ll 0$. Rephrased in a different way (using the fact that $Q$ is a bijection), we will show that 
\[0\neq Qx \;\underline{\gg}\; 0 \Rightarrow x\ll 0.\]
The condition $Qx \;\underline{\gg}\; 0$ gives us a system
\[ \begin{cases}
-v_1x_1+x_2\geq 0\\
x_1-v_2x_2+x_3\geq 0\\
\ldots\\
x_{n-1}-v_nx_n\geq 0
\end{cases}
\]
where all the $v_i$'s are $\geq 2$. Let $k$ be an index with 
\[x_k=\max_{i}\{x_i\}_i.\]
We want to show that $x_k<0$. Suppose that $1<k<n$. Then 
\begin{align*}
x_{k-1}-v_kx_k+x_{k+1}&\geq 0 \\
x_{k-1}+x_{k+1}&\geq v_kx_k\\
\end{align*}
and therefore
\[ 2x_k\overset{(a)}\geq x_{k-1}+x_{k+1} \geq v_kx_k\overset{(b)}\geq  2x_k.\]
The inequality $(a)$ follows by the definition of $x_k$, while $(b)$ is true if $x_k\geq 0$ (if it is $<0$ then we would be already done). This implies that $x_{k-1}=x_{k+1}=x_k$ and so we can assume, by iterating this argument, that $k=1$ (the case $k=n$ is the same).
We have:
\begin{align*}
-v_1x_1+x_2&\geq 0 \\
x_2&\geq v_1x_1\\
\end{align*}
and again, as before 
\[x_1\geq x_2\geq v_1x_1\geq 2x_1.\]
Therefore $x_1\geq 2x_1$, so $x_1\leq 0$. To exclude $x_1=0$ just notice that if this were the case, then from $-v_1x_1+x_2\geq 0$ it would follow that $x_2=0$ (being $x_1=0$ the maximum among the $x_i$'s), and consequently all the remaining $x_3=\ldots=x_n=0$, contradicting the assumption $Qx\neq 0$.
\end{prf}

\begin{prop}\label{concave}
For any rotation vector $x$ corresponding to a virtually overtwisted structure (i.e. with components $|x_i|\leq |y_i|$, with at least one strict inequality) we have
\[f(x)>f(y).\]
\end{prop}

\begin{prf}
Inside $\R^n$ we look at the region $D=\{(x_1,\ldots, x_n),\,|\,x_i|\leq |y_i|\,\, \forall i\}$. The goal is to show that the minimum of $\left.f\right|_D:D\to\R$ is realized on the vectors which correspond to the universally tight structures $y$ and $-y$, lying on $\partial D$.

Since $Q$ (and hence $Q^{-1}$) is negative definite, $f$ is concave. Being $f$ a negative definite norm, we know that it is has a unique maximum, which is the origin. Moreover, the minimum of $\left.f\right|_D$ is reached on the boundary $\partial D$. The fact that it is realized on $y$ and $-y$ follows from Lemma \ref{matrix}.
\end{prf}

This implies that the contact structures encoded by the vector $x$ cannot bound any Stein rational ball:

\begin{prf}[of Theorem \ref{minchi}]
Let $(X,J)$ be the Stein filling of $(L(p,q),\xi_{vo})$ described by the Legendrian realization of the chain of unknots associated with the vector of rotation numbers $x$.
By Proposition \ref{concave} we know that $f(x)>f(y)$, hence
\[c_1(X,J)^2=f(x)>f(y)=\sigma(X).\]
Since Equality \eqref{c2=sigma} is not satisfied, $(L(p,q),\xi_{vo})$ does not bound a Stein rational homology ball.
\end{prf}

\section{Coverings of tight structures on lens spaces and applications}\label{coveringsection}

In general, it can be hard to tell if the pullback of a tight contact structure on a 3-manifold along a given covering map is again tight. The situation is much easier if we restrict to lens spaces because of two reasons:
\begin{itemize}
\item[1)] tight structures are classified;
\item[2)] the fundamental group is finite cyclic, hence it is straightforward to determine their coverings.
\end{itemize}
If we start with a virtually overtwisted structure $\xi_{vo}$ on $L(p,q)$ we get an overtwisted structure $\pi^*\xi_{vo}$ on $S^3$, where $\pi:S^3\to L(p,q)$ is the universal cover. If $p$ is a prime number, then this is the only cover that $L(p,q)$ has, otherwise there is a bigger lattice of covering spaces depending on the divisors of $p$. 

Studying the behavior of coverings of a contact 3-manifold gives information about the fundamental group of its fillings:

\begin{thm} Let $Y$ be a closed and connected 3-manifold whose fundamental group $\pi_1(Y)$ is simple. Let $\xi$ be a virtually overtwisted contact structure on $Y$, and $(X,J)$ a Stein filling of $(Y,\xi)$. Then $X$ is simply connected.
\end{thm}

\begin{prf}
Let $i:Y\hookrightarrow X$ be the inclusion of the boundary $Y=\partial X$. Being $(X,J)$ a Stein filling of $(Y,\xi)$, the induced morphism $i_*:\pi_1(Y)\to \pi_1(X)$ is surjective. Moreover, by simplicity of $\pi_1(Y)$, we have that $\ker i_*$ can either be:
 
\begin{itemize} 
\item $\ker i_*=1:$ 

In this case, take a finite cover $p:(\widehat{Y},\xi_{ot})\to (Y,\xi)$ for which $\xi_{ot}$ is overtwisted. Call $n$ the degree of such cover. Define the group
\[G=i_*p_*\pi_1(\widehat{Y})\leq \pi_1(X),\]
consider the connected covering space of $X$ associated to $G$ and call it $\widehat{X_G}$. Since $\deg(\widehat{X_G}\to X)=n$, we have that $\widehat{X_G}$ is compact.
We are in the case where $i_*$ is an isomorphism and so $\partial \widehat{X_G}$ contains a diffeomorphic copy of $\widehat{Y}$. But by lifting the Stein structure from $X$ to $\widehat{X_G}$ we get a Stein structure on $\widehat{X_G}$ which fills the \emph{connected} contact boundary: note that any Stein semi-filling of a lens space is actually a filling, i.e. its boundary is connected (this comes from \cite[Theorem 1.4]{os}). Therefore we obtained a Stein filling of $\partial \widehat{X_G}=(\widehat{Y},\xi_{ot})$. This is not possible since the overtwisted contact structures are not fillable (as proved in \cite{eliashberg} and \cite{gromov}).

\item $\ker i_*=\pi_1(Y):$ 

This tells that $i_*$ is identically zero, and so that, by surjectivity, $\pi_1(X)=1$ as wanted.
\end{itemize} 
\end{prf}

\begin{cor}\label{pi1} Let $\xi$ be a virtually overtwisted structure on $L(p,q)$ with $p$ prime and let $(X,J)$ be one of its Stein fillings. Then $\pi_1(X)=1$.
\end{cor}

Now we want to study more carefully the behavior of the virtually overtwisted contact structures under covering maps, in order to derive some consequences on the possible fundamental groups of the fillings. The driving condition is the following observation:
let $p:\widehat{Y}\to Y$ be a covering map between compact and connected 3-manifolds, and let $i:Y\hookrightarrow X$ be the inclusion of the boundary $Y=\partial X$. Then, by covering theory:
\begin{align*}
\exists \mbox{ covering }\widehat{X}\to X \mbox{ that r}& \mbox{estricts to a covering } \partial\widehat{X}\to \widehat{Y} \\ 
&\Updownarrow \\
\ker i_* &\leq p_*\pi_1(\widehat{Y}).
\end{align*}
The way we want to apply this is to deduce that $\ker i_*$ should be big enough not to be contained in those subgroups of $\pi_1(Y)$ for which we can associate an \emph{overtwisted} cover. For example, if $X$ is a Stein filling of $Y$ and we are able to construct overtwisted coverings of $Y$ associated to every maximal subgroups of $\pi_1(Y)$, then the kernel of $i_*$ is forced to be the whole $\pi_1(Y)$, being this one the only subgroup of $\pi_1(Y)$ not contained in any maximal subgroups. 

By surjectivity of $i_*$ we would then conclude that $X$ is simply connected.

This looks to be promising because in the case of lens spaces it is easy to determine all the maximal subgroups of the fundamental group. It is nevertheless not so immediate to understand the behavior of the contact structure under the pullback map of a covering, but in certain cases we can use a necessary condition of compatibility of Euler classes to get some results. To better explain this, let us consider the following:

\begin{exmp} Let $(L(34,7),\xi_{vo})$ be obtained by contact $(-1)$-surgery on the Legendrian link of Figure \ref{hopf(5,7)}. If we orient the two components in the counter-clockwise direction we get rotation numbers respectively $+3$ and $-5$.

\begin{figure}[ht!]
\centering
\includegraphics[scale=0.5]{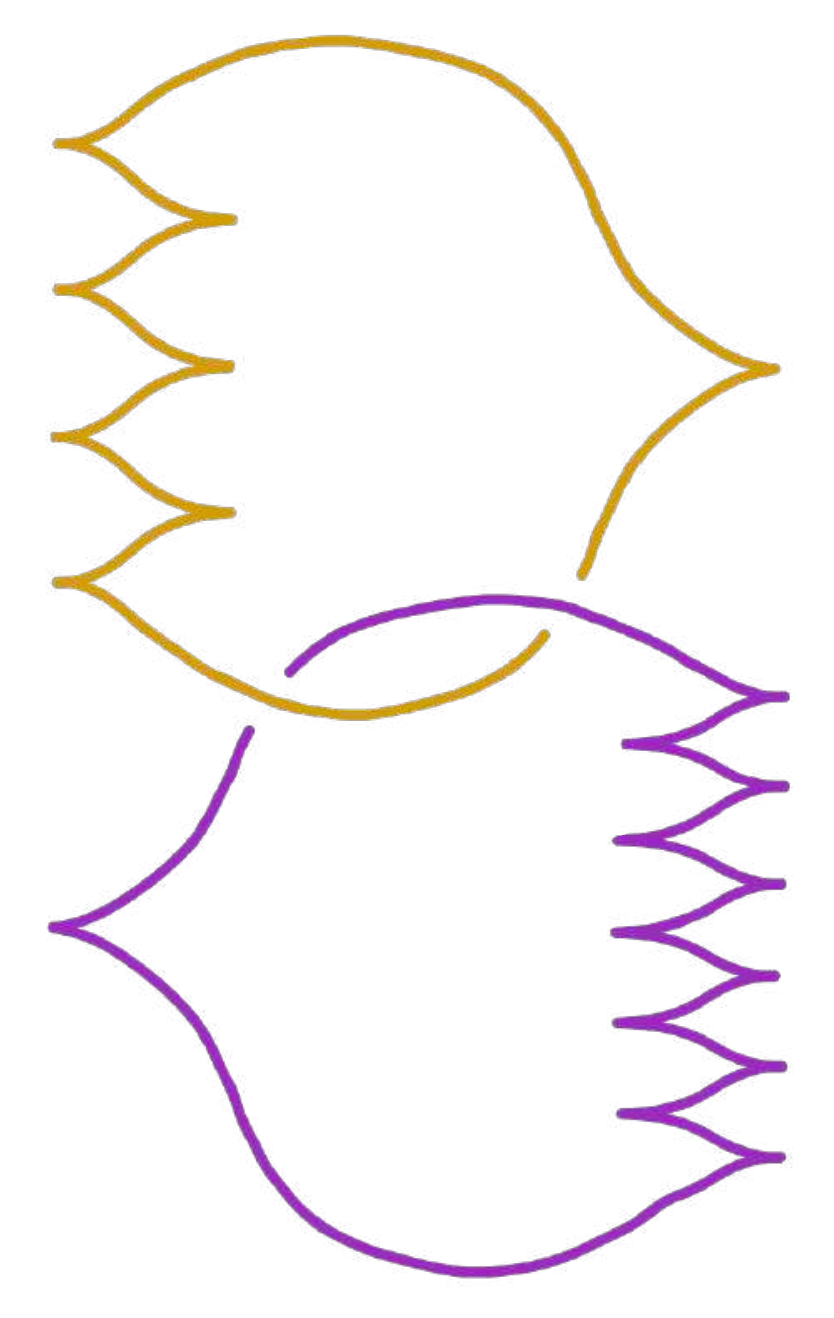}
\caption{Contact $(-1)$-surgery producing $L(34,7)$.}
\label{hopf(5,7)}
\end{figure}

After factoring $34=17\cdot 2$, we see that there are just two coverings:
\[L(17,7)\to L(34,7),\qquad L(2,7)\simeq L(2,1)\to L(34,7).\]
We will show that the given contact structure $\xi_{vo}$ on $L(34,7)$ lifts in both cases to an overtwisted structure. This tells us that, given any Stein filling $X$ of $L(34,7)$, the kernel on the inclusion map at the level of fundamental groups cannot be contained in $\Z/17\Z$ nor $\Z/2\Z$ and therefore is the whole $\Z/34\Z$, so $X$ is necessarily simply connected.
\begin{itemize}
\item The lift of $\xi_{vo}$ to $L(2,1)$ is overtwisted, because the only tight structure on $L(2,1)$ is universally tight and this one pulls backs to the tight structure on $S^3$, but since $\xi_{vo}$ is virtually overtwisted the lift to $S^3$ must be overtwisted.
\item To exclude that $\xi_{vo}$ pulls back to a tight structure on $L(17,7)$ we analyze the possible tight structures supported there. The fraction expansion of $17/7$ is
\[\frac{17}{7}=[3,2,4]\]
and so we see that there are 6 tight structures on $L(17,7)$ up to isotopy (and 3 up to contactomorphism, which are exhibited in Figure \ref{L(17,7)}).

For these structures we compute the Poincaré dual of the Euler class, viewed as an element of $\Z/17\Z\simeq H_1(L(17,7);\Z)$. The previous isomorphism is realized by choosing as a generator the meridian curve $\mu_1$ of the yellow curve with Thurston-Bennequin number $-2$.

Let $\xi$ be any of the three tight contact structures on $L(17,7)$ of Figure \ref{L(17,7)}. 
 
 \begin{figure}[h!]
\centering
\begin{subfigure}[t]{.3\textwidth}
  \centering
  \includegraphics[scale=0.4]{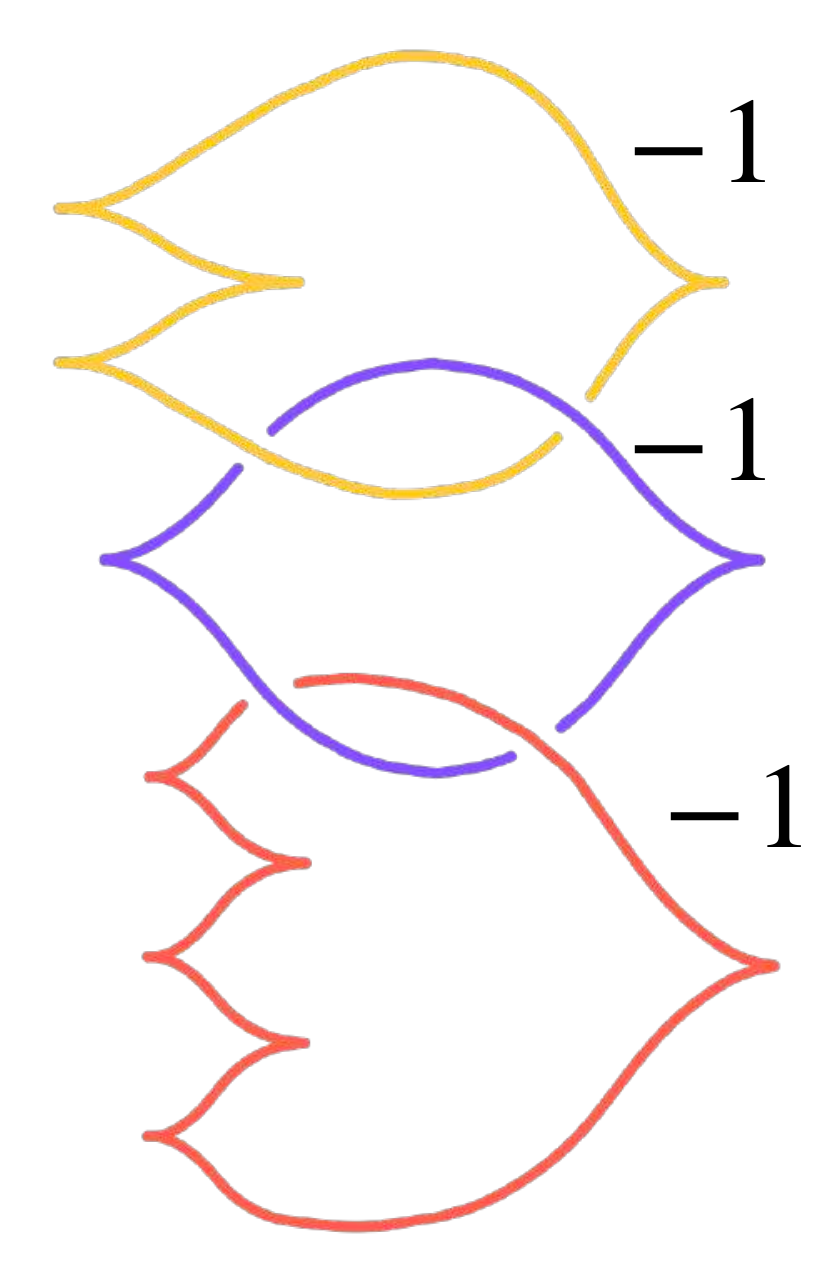}
  \caption{$\xi_1$}
  \label{xi1}
\end{subfigure}%
\begin{subfigure}[t]{.3\textwidth}
  \centering
 \includegraphics[scale=0.4]{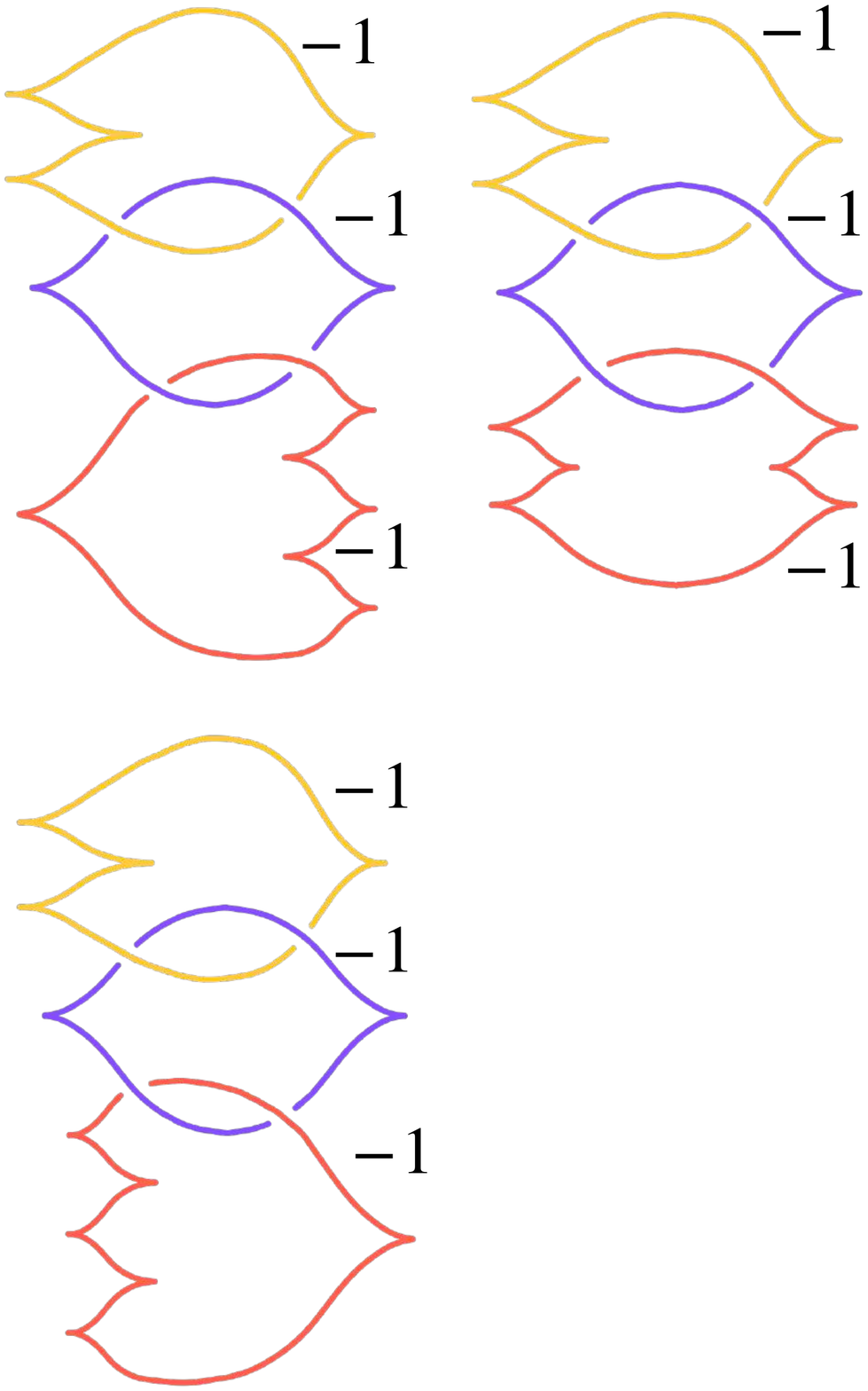}
  \caption{$\xi_2$}
  \label{xi2}
\end{subfigure}
\begin{subfigure}[t]{.3\textwidth}
  \centering
 \includegraphics[scale=0.4]{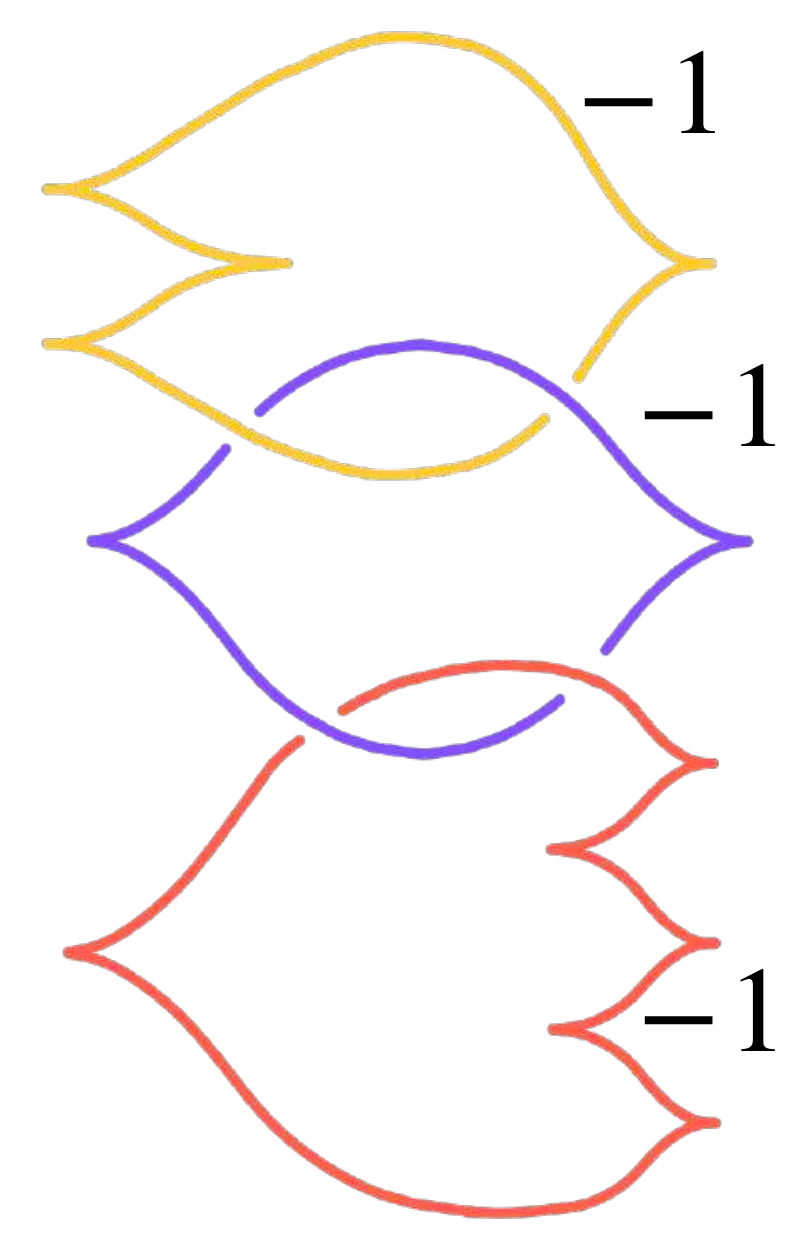}
  \caption{$\xi_3$}
  \label{xi3}
\end{subfigure}
\caption{Tight structures on $L(17,7)$.}
\label{L(17,7)}
\end{figure}

The class $\PD(e(\xi))$ is the image via the boundary map 
\[\partial: H_2(W,\partial W)\to H_1(\partial W)\simeq H_1(L(17,7))\]
of the Poincaré dual of the relative first Chern class of the Stein structure on $W$, where $W$ is the Stein domain described by the corresponding diagram of Figure \ref{L(17,7)}.
The Poincaré dual of the relative first Chern class is (see \cite[Proposition 8.2.4]{ozbagci})
\[\rot(K_1)[D_1,\partial D_1] +\rot(K_2)[D_2,\partial D_2] +\rot(K_3)[D_3,\partial D_3],\]
where the $K_i$'s are the three components of the link and the $[D_i,\partial D_i]$'s are the relative homology classes of the meridian disks of the 4-dimensional 2-handles attached to form the Stein filling $W$. Calling $\mu_i=\partial [D_i,\partial D_i]=[\partial D_i]$, for $i\in \{1,2,3\}$, the meridians of the attaching circles of these handles, we have

\[\PD(e(\xi))=\rot(K_1)\mu_1+\rot(K_2)\mu_2+\rot(K_3)\mu_3.\]

Let $Q$ be the matrix describing the intersection form of $W$, which is the same as the linking matrix
\[  \begin{bmatrix}
   -3 & 1 & 0 \\
 	1 & -2 & 1 \\
   0 & 1 & -4
   \end{bmatrix}.\]
From the exact sequence
\[\xymatrix{ H_2(W) \ar[r]^-Q & H_2(W,\partial W) \ar[r]^-{\partial} & H_1(\partial W)\simeq H_1(L(17,7))},\]
we get three linear relations
\[ \begin{cases}
-3	\mu_1+\mu_2=0\\
\mu_1-2\mu_2+\mu_3=0\\
\mu_2-4\mu_3=0.
\end{cases} \]
which tell us that $\mu_2=3\mu_1$ and $\mu_3=2\mu_2-\mu_1=5\mu_1$.
By putting everything together we get:

\begin{align*}
\PD(e(\xi))&=\partial(\PD(c_1(W,J)))\\
&=\partial (\rot(K_1)[D_1,\partial D_1] +\rot(K_2)[D_2,\partial D_2] +\rot(K_3)[D_3,\partial D_3])\\
&= \rot(K_1)\mu_1+\rot(K_2)\mu_2+\rot(K_3)\mu_3 \\
&=(\rot(K_1)+3\rot(K_2)+5\rot(K_3))\mu_1.
\end{align*}

If we substitute the values of the rotation numbers for the three different contact structures of Figure \ref{L(17,7)} we find:
\[ \PD(e(\xi_1))=11\mu_1, \qquad \PD(e(\xi_2))=\mu_1, \qquad \PD(e(\xi_3))=8\mu_1.\]

The contact structure described in Figure \ref{hopf(5,7)} we started from has 
\[ \PD(e(\xi))=12\mu, \]
with $\mu$ being the meridian curve of the yellow curve of Figure \ref{hopf(5,7)}.

Notice that $\mu$ is the image of the curve $\mu_1$ under the covering map $p:L(17,7)\to L(34,7)$. This is clear if we take the meridian curves of the single-component unknots with rational framing $-17/7$ and $-34/7$: here, the covering map is explicit on the two solid tori which form, once glued, each lens space, and the meridian upstairs is sent to the meridian downstairs. But when we expand from rational to integer surgery representation, we can think of keeping fixed the "first" solid torus, and glue in a series of thickened annuli, and at last the final solid torus. This is well described in \cite[Section 2.3]{saveliev}. This explains why $\mu_1$ is sent to $\mu$ by the covering map.

At the level of the homology group $H_1$ the covering map is a multiplication by 2 (the degree of the covering) and by naturality we need to find 
\[p_*(\PD(e(p^*(\xi))))=2\PD(e(\xi))\in H_1(L(34,7)).\]
But 
\[2\cdot 11\neq \pm 2\cdot 12 \in \Z/34\Z, \qquad 2\cdot 1\neq\pm  2 \cdot 12 \in \Z/34\Z, \qquad 2\cdot 8\neq\pm  2 \cdot 12 \in \Z/34\Z,\]
therefore we have that none of the three structures of Figure \ref{L(17,7)} is the pullback of our starting structure of Figure \ref{hopf(5,7)}. But those were the only (up to contactomorphism) tight structures on $L(17,7)$, so we conclude that the pullback is necessarily overtwisted, as wanted.

Note that we could have excluded a priori the contact structure $\xi_1$ of Figure \ref{xi1}, this being universally tight. 
\end{itemize}

Similar computations can be done if we start with a Legendrian representation of the Hopf link of Figure \ref{hopf(5,7)1} with rotation numbers $\pm(-3,1),\,\pm(-3,3),\,\pm(-3,5),\,\pm(-1,1),\,\pm(-1,3),\,\pm(-1,5)$. We made use of the software \emph{Mathematica} to carry out the computations and check that there is no tight structure on the double cover $L(17,7)$ with compatible Euler class.

\begin{figure}[h!]
\centering
\includegraphics[scale=0.6]{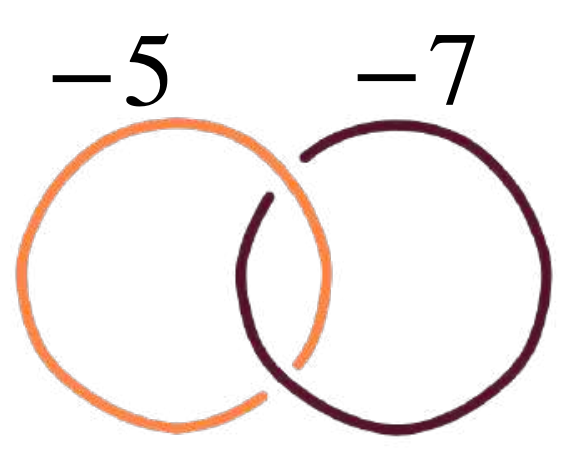}
\caption{Hopf link for $L(34,7)$.}
\label{hopf(5,7)1}
\end{figure}

The fact that the Stein fillings of these virtually overtwisted structure on $L(34,7)$ are simply connected can be deduced, as we just did, simply by looking at the two different coverings. This is something we already knew from the classification of fillings of those lens spaces obtained by contact surgery on the Hopf link, since the fraction expansion of $34/7$ has length 2 (see 	\cite{fossati}).
\end{exmp}

\begin{exmp}

Sometimes, an even quicker argument can be used to understand the behavior of a contact structure along certain covering maps. Let's take as an example $L(52,11)$, whose associated fraction expansion has length 3:
\[-\frac{52}{11}=[-5,-4,-3].\]
The two maximal subgroups of $\Z/52\Z$ are $\Z/4\Z$ and $\Z/26\Z$, and again, by running the computation of the Euler classes as above, we can determine which virtually overtwisted contact structure on the base cannot lift to a tight structure. 
But if we look at the covering of degree 13, we find $L(4,11)\simeq L(4,3)$ as total space, and since 
\[-\frac{4}{3}=[-2,-2,-2],\] 
we see that the only tight structure it supports is universally tight. 
Similarly, if we consider the covering $L(13,11)\to L(52,11)$ which has degree 4, we notice that
\[-\frac{13}{11}=[-2,-2,-2,-2,-2,-3]\] 
and hence also $L(13,11)$ supports only universally tight structures, among the tight ones. 
In the covering lattice of $L(52,11)$ it remains to study just the case of $L(26,11)$, for which the behavior can be more subtle (see next section, Theorem \ref{thmL5211}).

\[
\xymatrix{ & (L(13,11),{\color{red}\xi_{ot}})\ar[rr]^-{2:1}& & (L(26,11),{\color{blue}\xi_{?}})\ar[rd]_-{2:1} \\
(S^3,{\color{red}\xi_{ot}})\ar[ru]^-{13:1}\ar[rr]_-{2:1} & & (L(2,11),{\color{red}\xi_{ot}}) \ar[ru]_-{13:1} \ar[rd]_-{2:1}& & (L(52,11),{\color{blue}\xi_{vot}})\\
& & & (L(4,11),{\color{red}\xi_{ot}})\ar[ru]_-{13:1}
}
\]

\end{exmp}

\subsection*{A closer look to the coverings between lens spaces}

The test we made with the Poincaré duals gives only a necessary condition that does not guarantee that the pullback of a given tight contact structure is a tight contact structure simply because characteristic classes match.
So what can be said when there is compatibility between the Euler class of the contact structures of the base and of the covering?
We will try to present the idea of this subsection by starting from an example.

\begin{figure}[h!]
\centering
\begin{subfigure}[t]{.5\textwidth}
  \centering
  \includegraphics[scale=0.5]{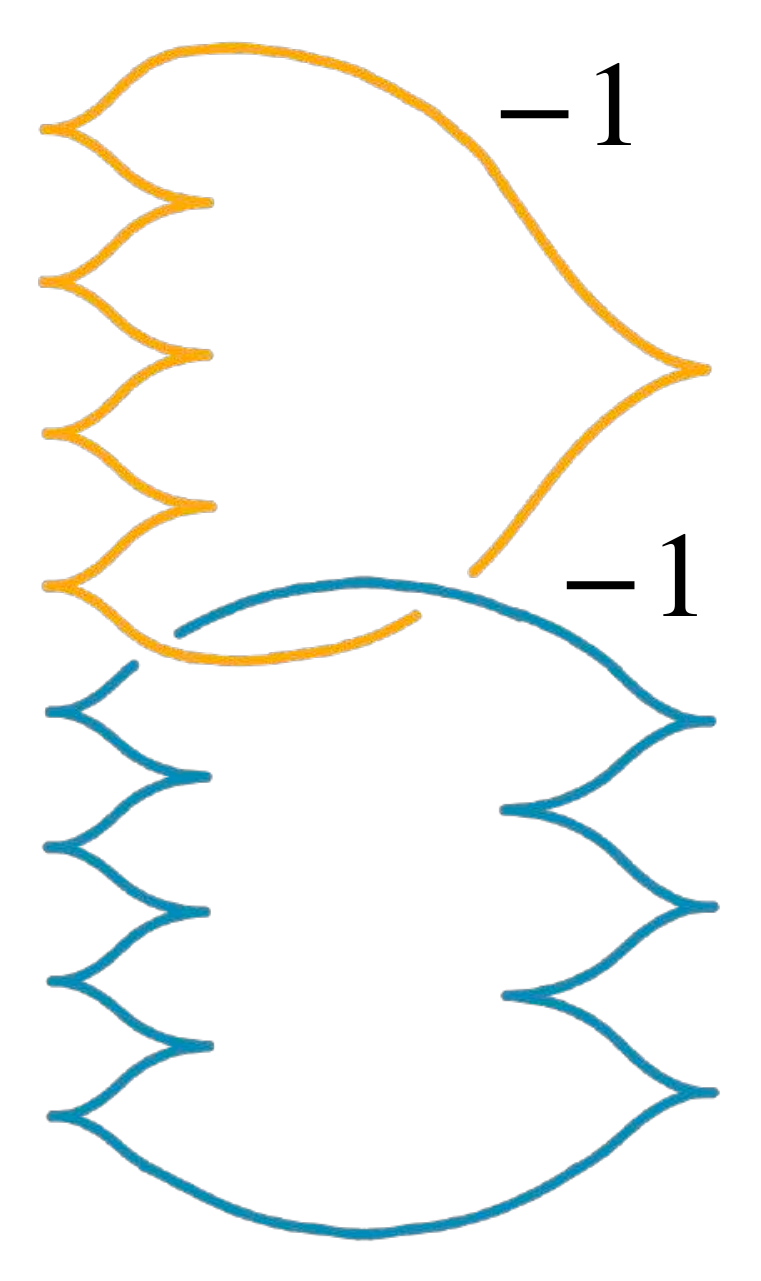}
  \caption{Contact structure on $L(34,7)$.}
  \label{L347}
\end{subfigure}%
\begin{subfigure}[t]{.5\textwidth}
  \centering
 \includegraphics[scale=0.5]{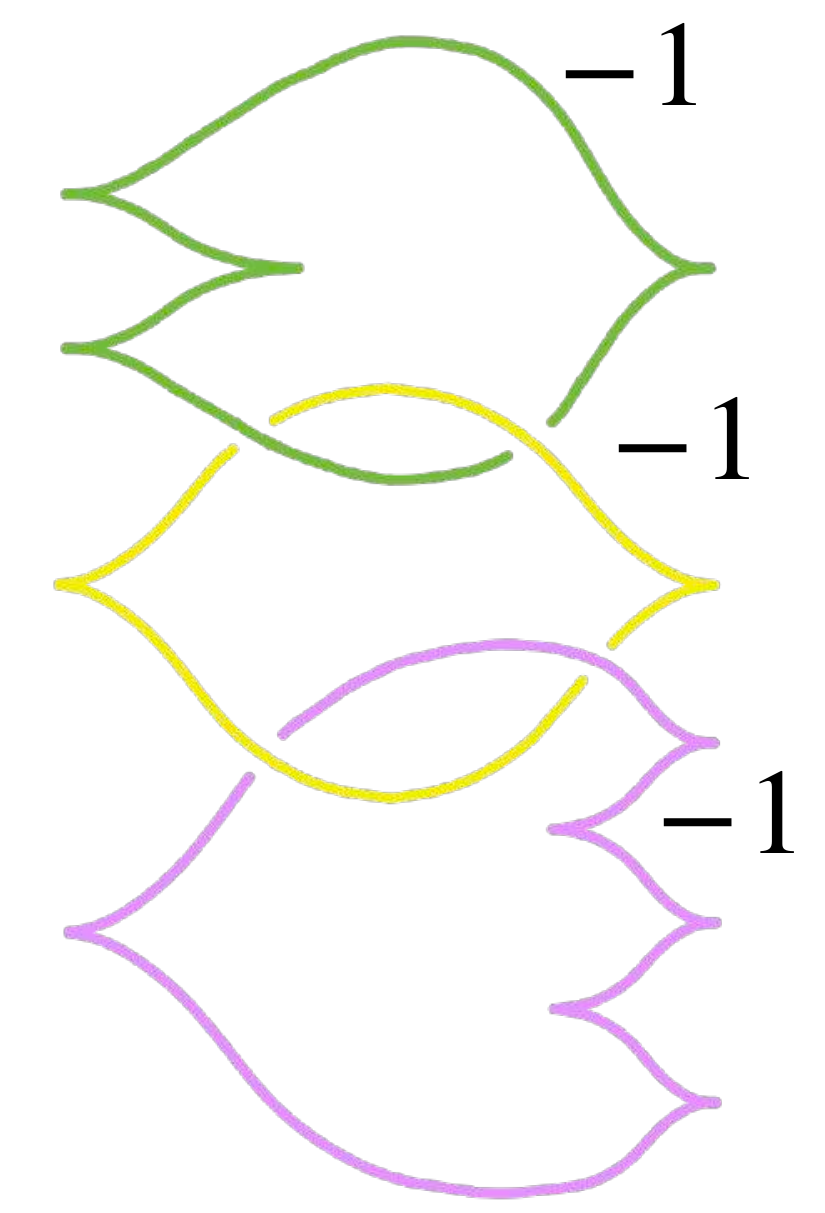}
  \caption{Contact structure on $L(17,7)$.}
  \label{L177}
\end{subfigure}
\caption{ }
\label{covering}
\end{figure}

Again, we choose to describe the double cover of $L(34,7)$. This time we fix the virtually overtwisted structure $\xi$ on $L(34,7)$ where the components of the link have rotation numbers $+3$ and $+1$ respectively, see Figure \ref{L347}. The computation shows that the Poincaré dual of the Euler class of $\xi$ is $+8\in \Z/34\Z$ (via the same identification of $H_1(L(34,7))\simeq \Z/34\Z$ as before). On the double cover $L(17,7)$ we take the tight structure $\hat{\xi}$ corresponding to the rotation vector $(1, 0, -2)$, as showed in Figure \ref{L177}. 

\noindent By running the computation, we find that $\PD(e(\hat{\xi}))=+8\in \Z/17\Z$, so that the covering map 
\[p:L(17,7)\to L(34,7)\]
takes $\PD(e(\hat{\xi}))$ to $16=2\PD(e(\xi))$, as it should certainly happen if $\hat{\xi}$ were isotopic to $p^*\xi$. But we will show that this is not the case, and argue that $p^*\xi$ is instead \emph{overtwisted}.

To do this, we need to use the description of tight structures on lens spaces of \cite{honda}, which we recall after the following definition:
\begin{defn} 
Given a contact 3-manifold $(Y,\xi)$, a \emph{contact vector field} $v$ on $Y$ is a vector field whose flow preserves the contact planes. A smooth surface $\S\subseteq (Y,\xi)$ is \emph{convex} if there exists a contact vector field $v$ on $Y$ transverse to $\S$.
The \emph{dividing set} of $v$ on $\S$ is defined as 
\[\Gamma=\{x\in\S\;\mid\;v(x)\in\xi_x\}.\]
\end{defn}

Giroux proved in \cite{girouxconv} that the dividing set is a 1-dimensional submanifold, whose isotopy type is independent of the choice of the contact vector field. We now focus on the case when $\S=T^2$. The diving set for a convex torus consists of an even number of parallel circles (see \cite[Section 5.1]{os}). By identifying $T^2$ with $\R^2/\Z^2$, we can talk about the \emph{slope} of these circles as a pair of numbers, which depends on the choice of the identification: when $T^2=\partial D^2\times S^1$, we use the meridian curve as one direction.

\paragraph{Honda's algorithm.}

\begin{figure}[h!]
\centering
\includegraphics[scale=0.5]{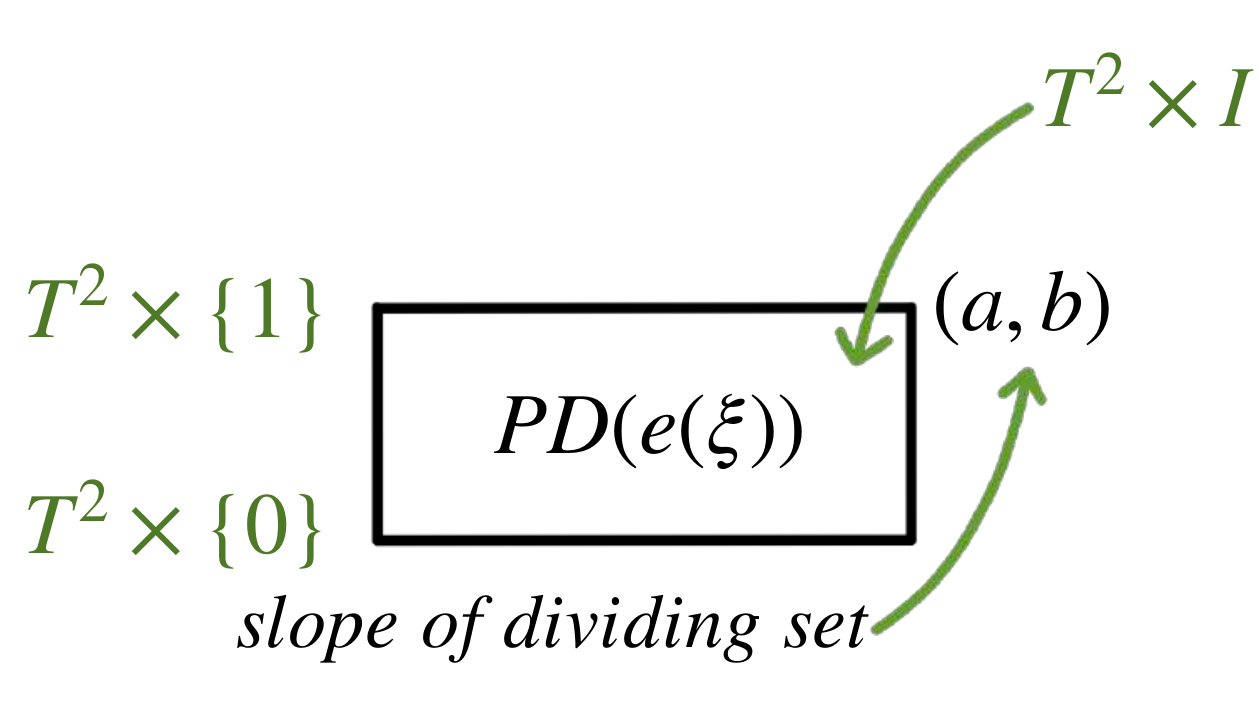}
\caption{Anatomy of a basic slice.}
\label{basicslice}
\end{figure}

In \cite[Section 4.3]{honda} it is explained how to cut a lens space, endowed with a tight contact structure, into two \emph{standard solid tori} and other pieces called \emph{basic slices}. With \emph{standard solid torus} we mean a small tubular neighborhood of a Legendrian knot, with standard coordinates on its boundary (see \cite[Section 2]{etnyrenotes}). On the other hand, a \emph{basic slice} is an oriented thickened torus $T^2\times I$ with a tight contact structure on it, such that the two boundary components are convex and satisfy certain conditions on the dividing sets. Each basic slice supports a unique tight contact structure, up to contactomorphism, but up to isotopy there are two classes: the isotopy class is determined by the sign of the (Poincaré dual of the) Euler class of the contact structure restricted to that basic slice. We always assume that the boundary tori are oriented according to the initial orientation on $T^2\times I$. A schematic picture of a basic slice is represented in Figure \ref{basicslice}.

\noindent The contact structure on the lens space is then encoded in the sequence of slopes on each basic slice and in the corresponding signs. We recall how the algorithm of Honda works for the lens space $L(p,q)$:

\noindent we start from the expansion $-p/q=[-a_1,\ldots,-a_n]$, with $a_i\geq 2$ for every $i$. Then we compute
\begin{align*}
-p_1/q_1:= & [-a_1,\ldots,-a_n+1] \\
-p_2/q_2:= & [-a_1,\ldots,-a_n+2] \\
-p_3/q_3:= & [-a_1,\ldots,-a_n+3] \\
\vdots
\end{align*}
until we get, after $k=a_n-1$ steps, to a rational number such that the length of its expansion is a number $m$, smaller than $n$, say
\[-p_k/q_k= [-b_1,\ldots,-b_m]. \]
This first set of numbers $\{-p_1/q_1,\ldots,-p_k/q_k\}$ will constitute the first block. Then we continue
\begin{align*}
-p_{k+1}/q_{k+1}:= & [-b_1,\ldots,-b_m+1] \\
-p_{k+2}/q_{k+2}:= & [-b_1,\ldots,-b_m+2] \\
\vdots
\end{align*}
until we get, after $h=b_m-1$ steps, to a rational number $-p_h/q_h$ such that the length of its expansion is less than $m$. This set of numbers $\{-p_k/q_k,\ldots,-p_h/q_h\}$ will constitute the second block. We go on this way until we reach the rational number $-1/1$. 

In total, we will produce an ordered set of blocks of ordered rational numbers which increase from $-p_1/q_1$ to $-1/1$, such that the numbers in each block have an associated continued fraction expansion of the same length. The boundary numbers, i.e. those which determine a change of length, appear twice: once at the bottom of a block, and then immediately after as the top of the following block. For example, the number $-14/3$ can appear twice, once as $[-5,-4,-1]$ and once as $[-5,-3]$. The expansion $[-5,-4,-1]$ determines the end of the block with length 3, while $[-5,-3]$ determines the start of the block of length 2. We record these rational numbers $-p_i/q_i$ as pairs $(-q_i,p_i)$, simplified to the form where $q_i$ is as small as possible. These numbers correspond to the slope of the dividing sets of the contact structure under analysis, when restricted to the corresponding basic slice.

Then we remove a standard torus from the lens space $(L(p,q),\xi)$ and we picture what is left in the following way: we draw the basic slices starting from the slope $-p_1/q_1$ until $-1/1$, divided into the blocks as described above. At the end of this thickened torus we draw the other basic torus. 

As we explained, every basic slice comes equipped with boundary slopes described by two rational numbers, which are represented by pairs $(-q,p)$ and $(-q',p')$. Honda proved that taking the difference of these values gives the Poincaré dual of the Euler class restricted there, up to sign, as an element of $H_1(T^2)\simeq \Z\oplus\Z$, written in the basis $(\partial D^2,S^1)$ specified by the lower solid torus. As mentioned above, the isotopy class of the unique (up to contactomorphism) contact structure on each basic slice is specified by the sign of the restriction of the Poincaré dual of the Euler class. Within a single block of basic slices, the only thing that matters is how many positive and negative signs we have, but not where these are placed: this is a consequence of a property of shuffling, which says that rearranging the signs within a block gives an isotopic contact structure, see \cite[Section 4.4.5]{honda}. This is coherent with the fact that, when drawing a Legendrian unknot with its stabilizations, we don't need to remember if we first stabilized positively or negatively, but just the final result. 

\textbf{To sum up}: if we start from a chain of $n$ Legendrian unknots, we get $n$ blocks (one for every component) of $a_i-2$ basic slices each (where $-p/q=[-a_1,\ldots ,-a_n]$). Every positive/negative stabilization that we see in the Legendrian link corresponds to a plus/minus in the corresponding block. Notice that when a coefficient in the expansion is $-2$, then its corresponding block will be empty, reflecting the fact that there is no choice of placing stabilizations in a Legendrian knot with Thurston-Bennequin number $-1$.

\vspace{0.5cm}

\begin{figure}[h!]
\centering
\begin{subfigure}[t]{.5\textwidth}
  \centering
  \includegraphics[scale=0.5]{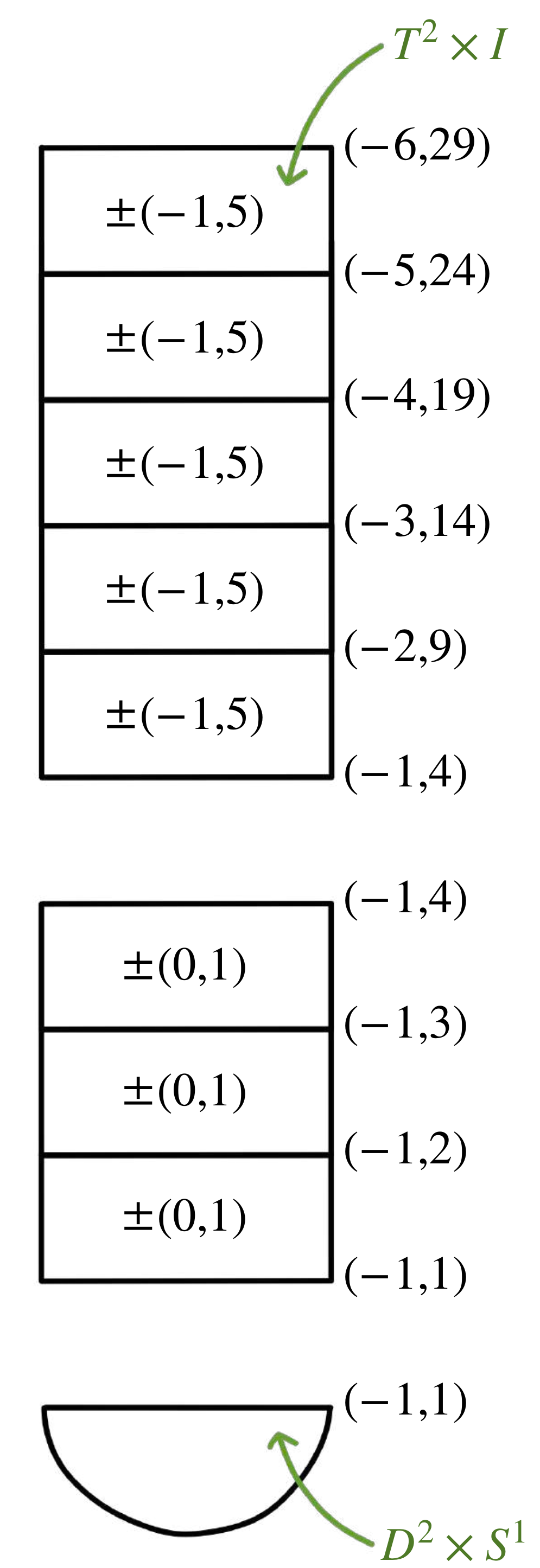}
  \caption{Subdivision into basic slices.}
  \label{blocks}
\end{subfigure}%
\begin{subfigure}[t]{.5\textwidth}
  \centering
 \includegraphics[scale=0.5]{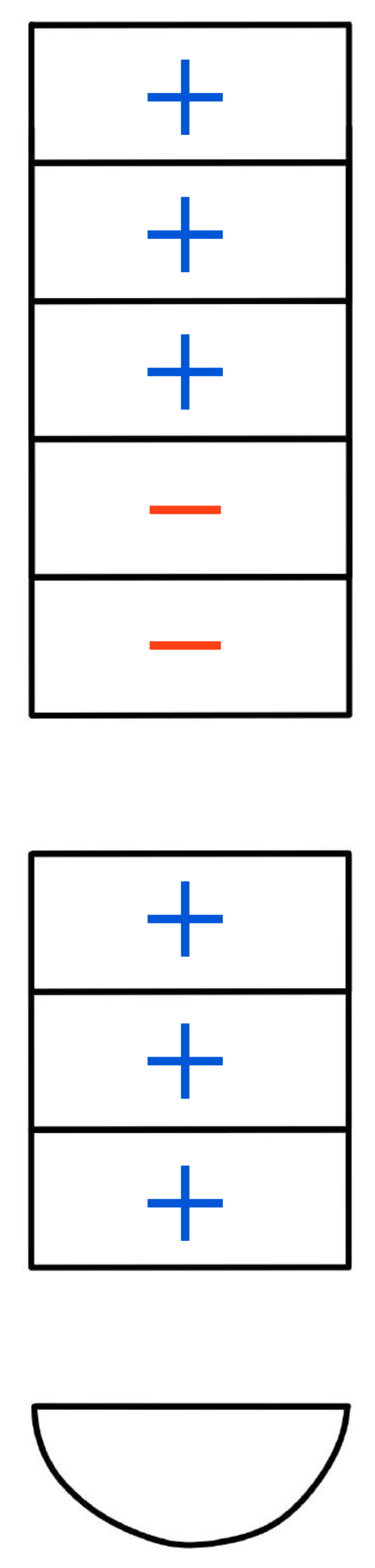}
  \caption{Signs of the basic slices.}
  \label{signs}
\end{subfigure}
\caption{Description of the contact structure $\xi$ on $L(34,7)$.}
\label{blockssigns}
\end{figure}

For example, see Figure \ref{blocks}: the algorithm gives two blocks of 5 and 3 basic slices respectively, where the slopes of the dividing sets on the boundary are indicated there.

From this picture it is also easy to calculate the Poincaré dual of the Euler class of the structure we choose according to the signs of each basic slice (indicated with colors blue and red in Figure \ref{signs}). By capping off with the solid torus below, we make the first $S^1$-factor of $T^2\times I$ nullhomologous, so we can just focus on the second entry in homology. The Poincaré dual of the Euler class of the structure is finally understood in the first homology group of the lens space once we glue the other solid torus (above). The structure in Figure \ref{blockssigns} has $\PD(e)$ given by 
\begin{equation}\label{PDdown}
5+5+5-5-5+1+1+1=8\in \Z/34\Z,
\end{equation}
and it is exactly the one resulting from contact $(-1)$-surgery on the Legendrian Hopf link of Figure \ref{L347}, where the component with Thurston-Bennequin number $-4$ has rotation $+3$ (corresponding to the three pluses in the lower block), and the other one has rotation $+1$ (corresponding to the upper block with three pluses and two minuses).

\begin{figure}[h!]
\centering
\begin{subfigure}[t]{.5\textwidth}
  \centering
  \includegraphics[scale=0.5]{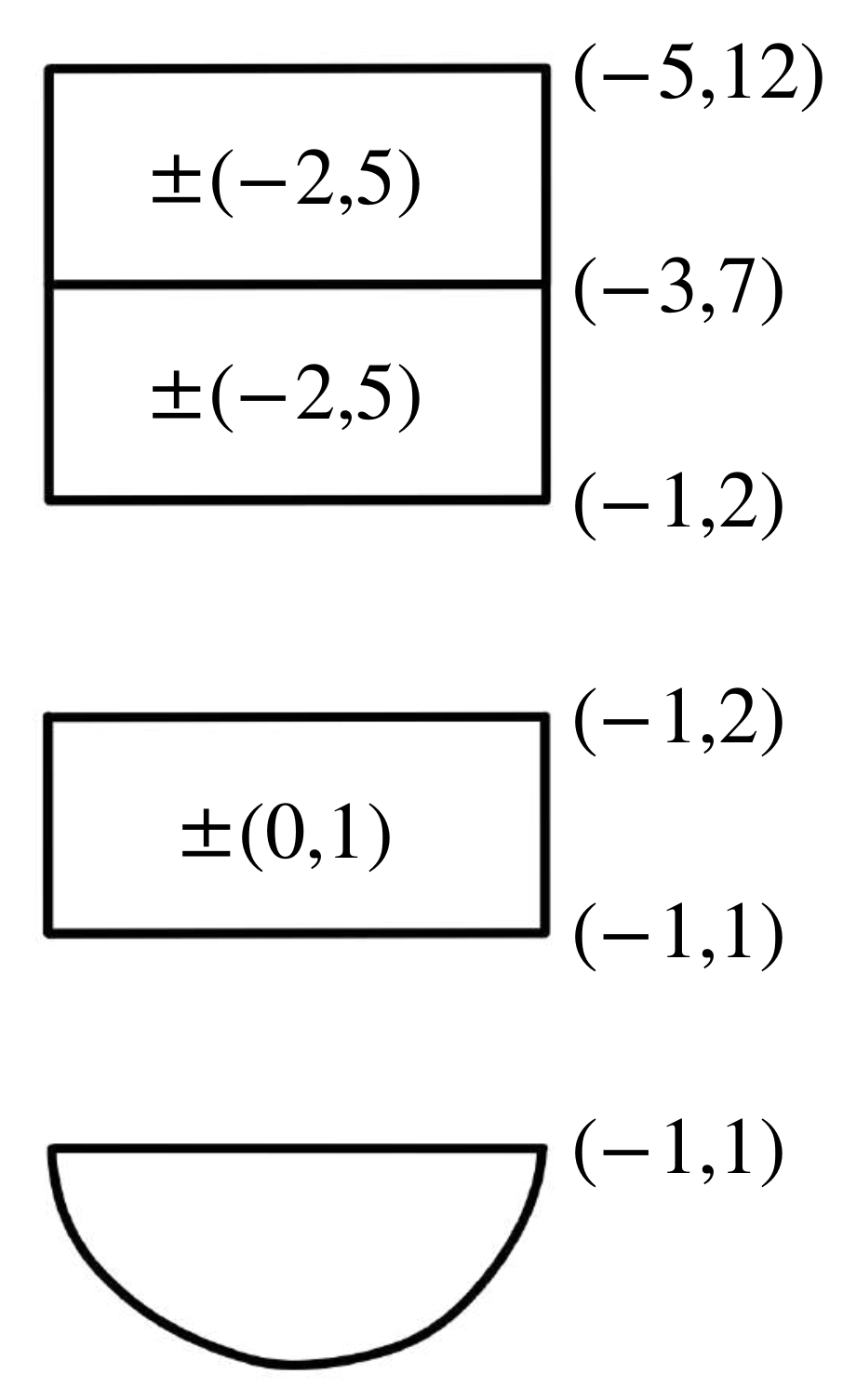}
  \caption{Subdivision into basic slices.}
  \label{blocks2}
\end{subfigure}%
\begin{subfigure}[t]{.5\textwidth}
  \centering
 \includegraphics[scale=0.5]{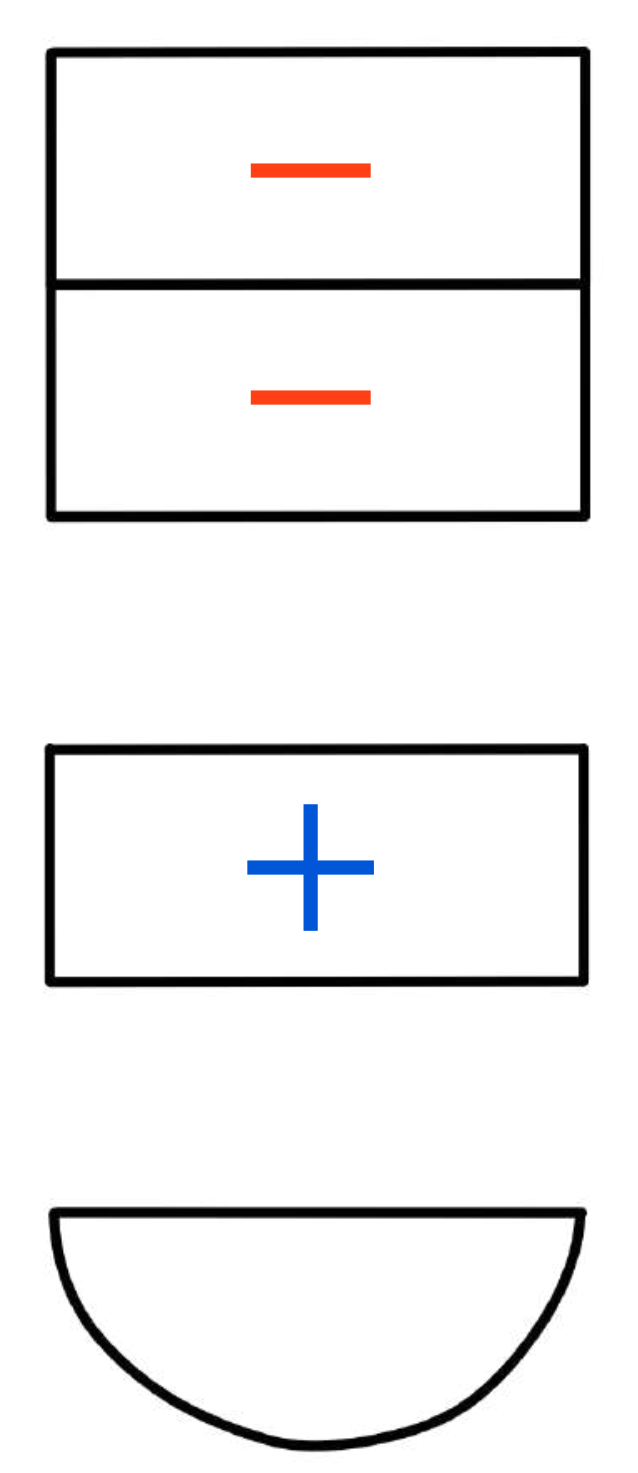}
  \caption{Signs of the basic slices.}
  \label{signs2}
\end{subfigure}
\caption{Description of a contact structure on $L(17,7)$.}
\label{blockssigns2}
\end{figure}

Now we look at the double covering map, which, on every basic slice, looks like
\[(z,w)\to (z,w^2),\]
where $z$ is the coordinate corresponding to $\partial D^2$ (which will be capped off when the lower solid torus is glued), and $w$ is the coordinate of the other $S^1$-factor.
Honda asserts that the contact structure on $L(p,q)$ can be encoded into one solid torus (the other one being standard), which is the one we chose to draw. This is why to analyze the behavior of the contact structure along the covering map we need to understand how this solid torus lifts. So we split $L(17,7)$ with a tight structure into two solid tori: the first one is pictured in Figure \ref{blocks2}, and subdivided into a block of two basic slices, plus a single basic slice, plus a standard solid torus; the other solid torus is a standard torus which will be glued on top of the uppermost basic slice, and which is not pictured. 

All the tight contact structures on $L(17,7)$ are encoded in the decomposition of the represented solid torus into these pieces: by choosing the sign of the basic slices we produce all the different (up to isotopy) 6 tight structures that $L(17,7)$ supports. As a double check, one can think at the different Legendrian representatives of the 3-components link made by a chain of unknots with Thurston-Bennequin numbers $-2$, $-1$ and $-3$ (notice indeed that $-17/7=[-3,-2,-4]$). The candidate tight contact structure on $L(17,7)$ which should be the pullback of the one on $L(34,7)$ described by Figure \ref{blocks} has the single basic slice with positive sign, and the other two in the block with negative signs. This corresponds to the choice of the rotation numbers for the components of the link to be $+1$, $0$ and $-2$: the link on which contact $(-1)$-surgery should give the pullback structure on $L(17,7)$ along the covering map is pictured in Figure \ref{L177}. The reason why this is the correct candidate is because this is the only case where we have compatibility of Euler classes: the computation (which can be performed in two different ways) shows that the Poincaré dual of the Euler class upstairs is $-9\equiv 8 \pmod{17}$, which gets sent to $16=2\cdot 8 \in \Z/34\Z$, which, as we already computed in Equation \eqref{PDdown}, is the double of the Poincaré dual of the Euler class downstairs.

\begin{figure}[h!]
\centering
\begin{subfigure}[t]{.5\textwidth}
  \centering
  \includegraphics[scale=0.5]{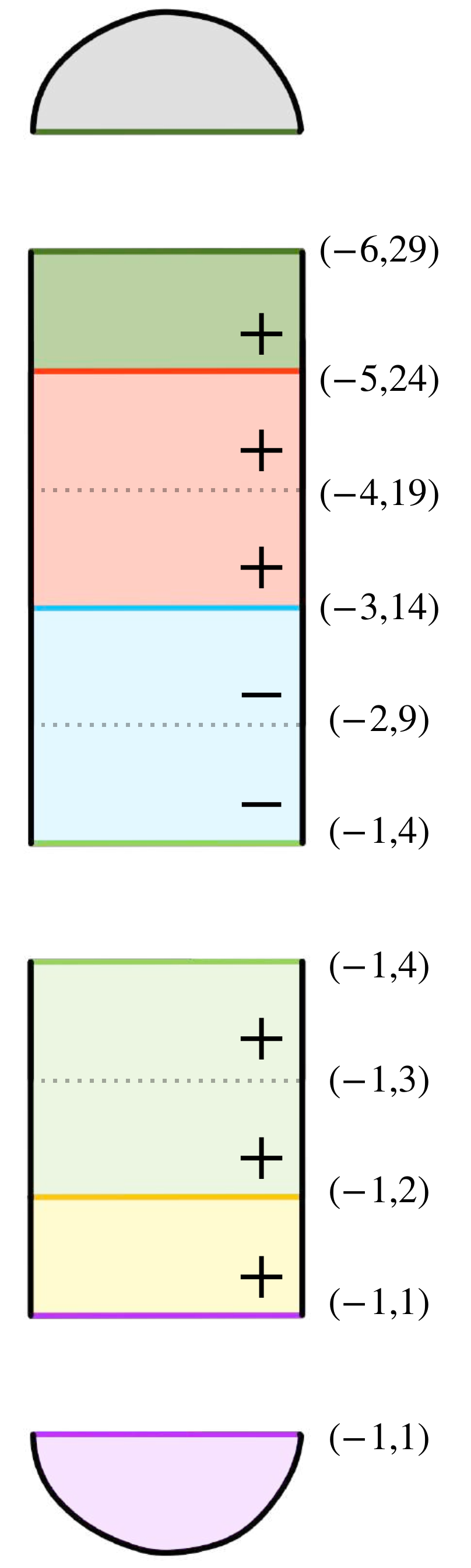}
  \caption{$L(34,7)$}
  \label{colors347}
\end{subfigure}%
\begin{subfigure}[t]{.5\textwidth}
  \centering
 \includegraphics[scale=0.5]{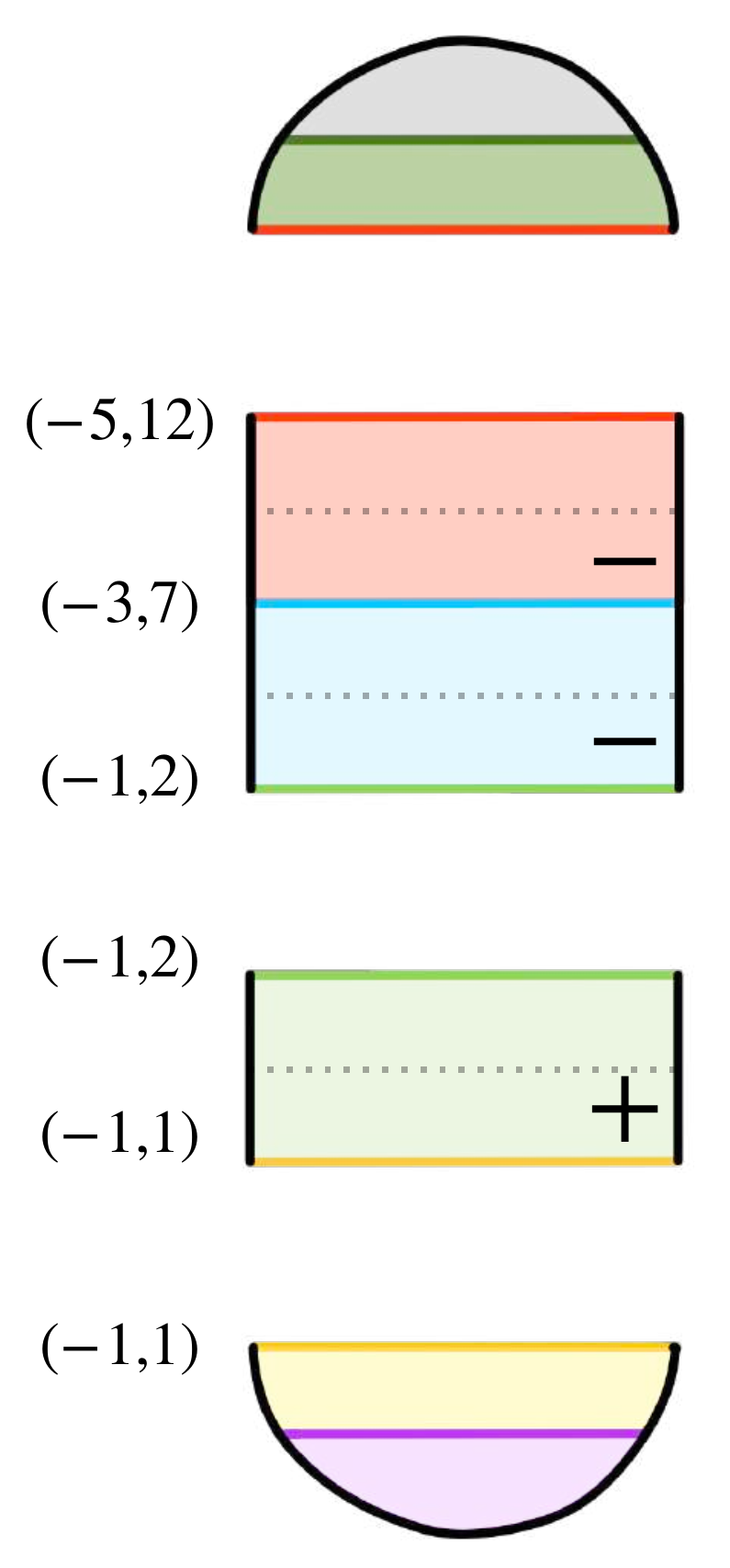}
  \caption{$L(17,7)$}
  \label{colors177}
\end{subfigure}
\caption{Behavior of slices under the covering map.}
\label{colors}
\end{figure}

But now we argue that there cannot be compatibility in the signs of the basic slices of $L(34,7)$ and $L(17,7)$. Indeed, once a sign for a basic slice downstairs is chosen, then its lift should have the same sign (see \cite[Section 1.1.4]{honda2}).
By lifting the dividing sets of the various convex tori we see where the different basic slices go: Figure \ref{colors} is describing this by means of colors. Notice that the lowest basic slice of $L(34,7)$ is pulled back inside the standard torus, and the same is true for the uppermost slice. Therefore the behavior of the contact structure upstairs is regulated by what happens to the central slices, i.e. from the yellow line $(-1,2)$ to the red line $(-5,24)$.

But here we finally see the contradiction. While:
\begin{itemize}
\item[1)] the positive slices from yellow $(-1,2)$ to green $(-1,4)$ lift to a positive slice in $L(17,7)$ and
\item[2)] the negative slices from green $(-1,4)$ to blue $(-3,14)$ lift to a negative slice in $L(17,7)$, we have that
\item[3)]  the \emph{positive} slices from blue $(-3,14)$ to red $(-5,24)$ lift to a \emph{negative} slice in $L(17,7)$
\end{itemize}
and this is not possible. No matter how we decide to shuffle the basic slices in each single block (see \cite[Section 4.4.5]{honda}), we always end up with a contradicting situation (as proved in Theorem \ref{thmcoverovertwisted}).

This tells us that, even if there is a tight virtually overtwisted structure on $L(17,7)$ whose Euler class is compatible with the  structure $\xi$ we chose on $L(34,7)$, the pullback of $\xi$ along the double covering map is \emph{overtwisted}, as claimed.

\begin{thm}\label{thmcoverovertwisted} Any virtually overtwisted structure on $L(34,7)$ lifts to an overtwisted one along the double cover 
\[L(17,7)\to L(34,7).\]
\end{thm}

\begin{prf}
We argue here using the behavior of the basic slices described in Figure \ref{colors}.
Look at the three basic slices in $L(17,7)$, Figure \ref{colors177}, regardless of the signs. Call $\hat{\xi}$ the pullback of a given $\xi$ on $L(34,7)$ and compare the Poincaré dual of their Euler classes. Assuming that the structures are both tight, we see that the choice of the sign of the red basic slice in $L(17,7)$ contributes to a $\pm 5$ for $\PD(\hat{\xi})$ and, pushed down, to a $\pm 10$ for $\PD(\xi)$. The same is true for the light blue slice, while the green slice gives a $\pm 1$ for $\PD(\hat{\xi})$ and a $\pm 2$ for $\PD(\xi)$. Moreover, inside $L(34,7)$ we have two extra slices (dark green and yellow in Figure \ref{colors347}), whose signs can be chosen independently. Requiring compatibility of Euler classes means to impose
\[\PD(\xi)\equiv \PD(\hat{\xi})\;\pmod{17}.\]
Therefore, according to what we have just said:
\[{\color{red}\pm 10}\, {\color{cyan}\pm 10} \, {\color{green}\pm 2}\,  {\color{olive}\pm 5}\, {\color{orange}\pm 1}\equiv {\color{red}\pm 5}\,	{\color{cyan}\pm 5} \, {\color{green}\pm 1} \;\pmod{17}\]
which is the same as
\[{\color{red}\pm 5}\, {\color{cyan}\pm 5} \, {\color{green}\pm 1}\,  {\color{olive}\pm 5}\, {\color{orange}\pm 1}\equiv 0\;\pmod{17}.\]
Clearly, this can be done only in two ways, namely by choosing all pluses or all minuses. And these correspond exactly to the two universally tight structures, for which we already knew that there is compatibility. Therefore, among the virtually overtwisted structures there cannot be a coherent choice of signs resulting in compatible Euler classes.
\end{prf}

\begin{thm}\label{thmL5211} Any virtually overtwisted structure on $L(52,11)$ lifts to an overtwisted one along all of its non-trivial covers.
\end{thm}

\begin{prf} At the end of previous section we argued that in the covering lattice of $L(52,11)$ the only case which was more subtle to describe was the double cover
\[L(26,11)\to L(52,11),\]
because otherwise we already knew that virtually overtwisted structures on the base would lift to overtwisted structures. We analyze this remaining case as we did before, by looking for compatibility between the signs of the basic slices and the count of the possible Euler classes.
Figure \ref{L5211L2611} shows where the basic slices go, from $L(26,11)$ to $L(52,11)$. 
\begin{figure}[h!]
\centering
\begin{subfigure}[t]{.5\textwidth}
  \centering
  \includegraphics[scale=0.5]{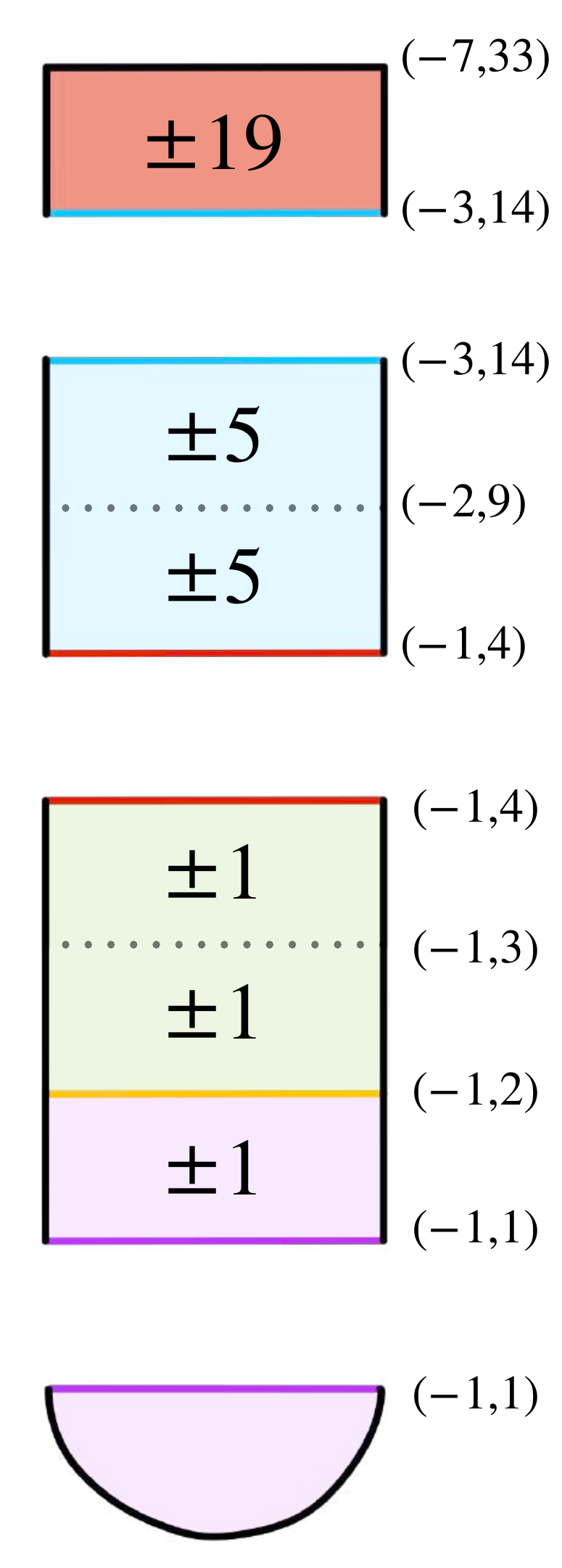}
  \caption{$L(52,11)$}
  \label{blocksL5211}
\end{subfigure}%
\begin{subfigure}[t]{.5\textwidth}
  \centering
 \includegraphics[scale=0.5]{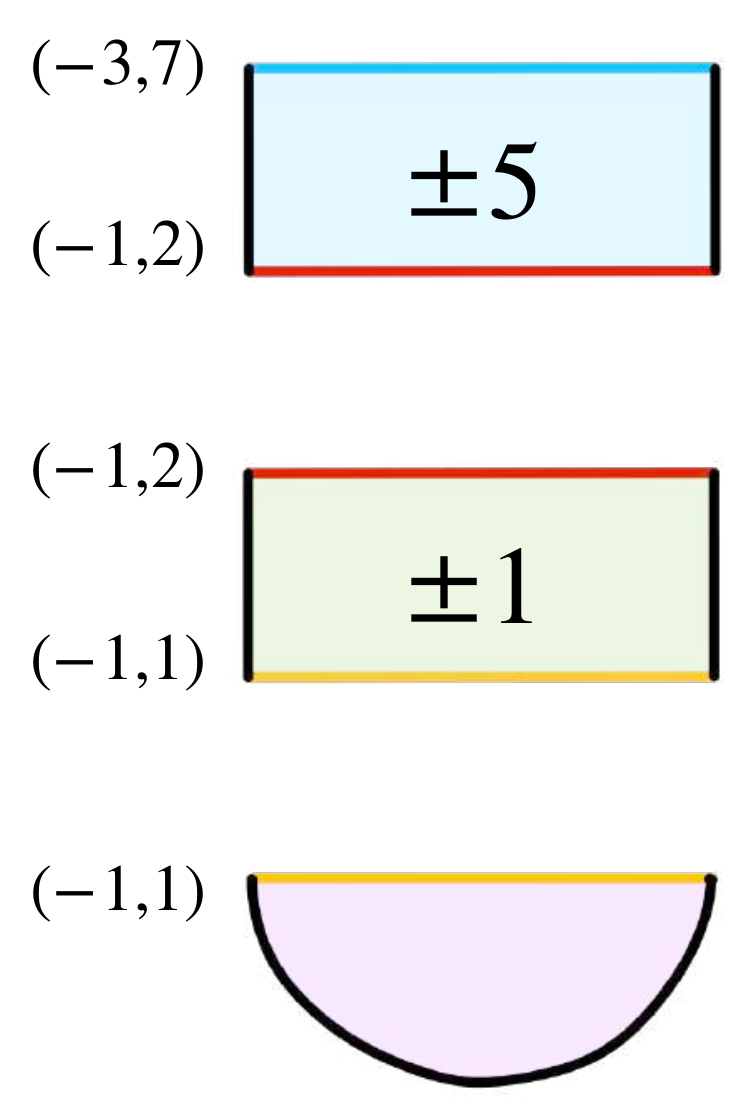}
  \caption{$L(26,11)$}
  \label{blocksL2611}
\end{subfigure}
\caption{Behavior of slices under the covering map.}
\label{L5211L2611}
\end{figure}
The count of the Poincaré duals of the two Euler classes gives
\[{\color{cyan}\pm 10}\, {\color{green}\pm 2} \, {\color{red}\pm 19}\,  {\color{violet}\pm 1} \equiv {\color{cyan}\pm 5}\,	{\color{green}\pm 1} \;\pmod{26}\]
which is the same as
\[{\color{cyan}\pm 5}\, {\color{green}\pm 1} \, {\color{red}\pm 19}\,  {\color{violet}\pm 1} \equiv 0 \;\pmod{26}.\]
Again, we see that this can be done only in two ways, namely by choosing all pluses or all minuses, which correspond exactly to the two universally tight structures. Therefore, among the virtually overtwisted structures there cannot be a coherent choice of signs resulting in compatible Euler classes.
\end{prf}

We can finally give a proof of Theorem \ref{thmcovering}, which stated that every virtually overtwisted contact structure on $L(p,q)$ lifts along a degree $d$ covering to a structure which is overtwisted, if we are in the condition where $q<p<dq$. 

\begin{prf}[of Theorem \ref{thmcovering}]
If the pullback of the contact structure were tight, it should fit with the description of tight structures according to the basic slices subdivision. We claim that the lower solid torus $\mathcal{H}_1$ until the level $-p/q$ gets all pulled back into the standard solid torus whose dividing set has slope $-1/1$. This comes from the fact that the curve with slope $(-q,p)$ pulls back to the one with slope $(-dq,p)$, according to the behavior
\[\xymatrix{
-\frac{p}{dq}\ar@{|->}[r]^-{\cdot d}&-\frac{p}{q}. }\]
But by assumption $-\frac{1}{1}<-\frac{p}{dq}$, and since the slopes of the dividing sets are increasing when read from top to bottom (compare with Figure \ref{blocks}), the claim follows.

\noindent Since we are considering a virtually overtwisted structure on $L(p,q)$, the pullback of $\mathcal{H}_1$ \emph{cannot be tight}, otherwise it would be universally tight, being it a subset of a solid torus in standard coordinates (which does not support virtually overtwisted structures). Therefore, we must have here an overtwisted disk, as wanted.
\end{prf}

\begin{cor} Let $p_1,\,p_2$ be prime numbers, not necessarily distinct, and let $q$ be an integer such that $p_i<q<p_1p_2$ for $i=1,\,2$. Then each non-trivial covering of $(L(p_1p_2,q),\xi_{vo})$ is overtwisted, for any virtually overtwisted structure $\xi_{vo}$.
\end{cor}

\begin{prf}
It is a direct consequence of previous theorem, since 
\[p_1p_2<dq,\]
where $d$ is either $p_1$ or $p_2$ (which are the only possible degrees for a non-trivial covering).
\end{prf}

\begin{nrem} The hypothesis of Theorem \ref{thmcovering} can be relaxed by just requiring that $p'<dq'$, where $p'$ and $q'$ are determined as follows: let 
\[-\frac{p}{q}=[-a_1,\ldots ,-a_n]\]
be the continuous fraction expansion, with $a_i\geq 2$ for each $i=1,\ldots, n$. Then define $p'$ and $q'$ as
\[[-a_1,\ldots ,-a_n+1]=-\frac{p'}{q'}.\]
In this way we have
\[-\frac{p}{q}<-\frac{p'}{q'}\]
so that the requirement $-1<-p'/q'$ is less restrictive. The reason why Theorem \ref{thmcovering} stays true with this weaker assumption is that the description of a contact structure via basic slices shows as the smallest slope (hence on top of the uppermost block) precisely the slope $-p'/q'$ (see \cite[Section 4.6]{honda}). To ask that, from this level down, the solid torus is pulled back inside the standard torus in the covering guarantees the existence of an overtwisted disk in the covering space, as argued in the proof of Theorem \ref{thmcovering}. 

There is another description of the two numbers $p'$ and $q'$ which is intrinsic in the sense that does not involve the computation of the continued fraction expansion: given $p$ and $q$, let $q^*$ be the multiplicative inverse of $q$, modulo $p$, i.e. $0<q^*<p$ and
\[q^*q\equiv 1\pmod{p}.\]
If we put $p'=p+q^*$, then $q'$ is the multiplicative inverse of $q^*$, modulo $p'$, i.e.
\[q'q^*\equiv 1\pmod{p'}.\]
\end{nrem}

\subsection*{Comparing $\pi_1$ and $\chi$ of a filling}

The goal of this section is to see some applications to concrete examples of Theorem \ref{pi1b2}, which is proved below.

\begin{prf}[of Theorem \ref{pi1b2}] 
Take the universal covering $\widetilde{X}\to X$, of degree $d$, whose boundary is the (connected) covering $L(p',q')\to L(p,q)$ of  degree $d$.
The Euler characteristics satisfy
\[\chi(\widetilde{X})=d\chi(X)\]
and hence, by Theorem \ref{maxchi},
\[\chi(X)=\frac{\chi(\widetilde{X})}{d}\leq \frac{1+l'}{d}.\]
\end{prf}

\begin{cor} Let $X$ be a Stein filling of a lens space $L(p,q)$ with a virtually overtwisted structure, and let $d$ be a divisor of $p$. If
\[2d>1+\length((p/d)/q),\]
then the fundamental group of $X$ cannot be $\Z/d\Z$.
\end{cor}

\begin{prf}
It follows by contradiction from Theorem \ref{pi1b2} if we look at the associated $d$-covering $\widehat{X}\to X$ and remember that $2\leq \chi(X)$, as proved in Section \ref{lowerboundsection} and \cite[Proposition A.1]{gollastar}. The number $\length((p/d)/q)$ has to computed after reducing $q$ modulo $p/d$.
\end{prf}

Sometimes, depending on the arithmetic of the rational numbers, it happens that the behavior of the basic slices of a covering is never compatible with the choice of signs determining the Euler classes, and this guarantees the covering itself to be overtwisted, which in turn implies that all the fillings are simply connected. But there are cases where a non-trivial cover of a tight virtually overtwisted structure stays as such, and so we need other arguments to calculate the fundamental group of a filling.

\begin{figure}[ht!]
\begin{subfigure}[t]{.5\textwidth}
  \centering
  \includegraphics[scale=0.5]{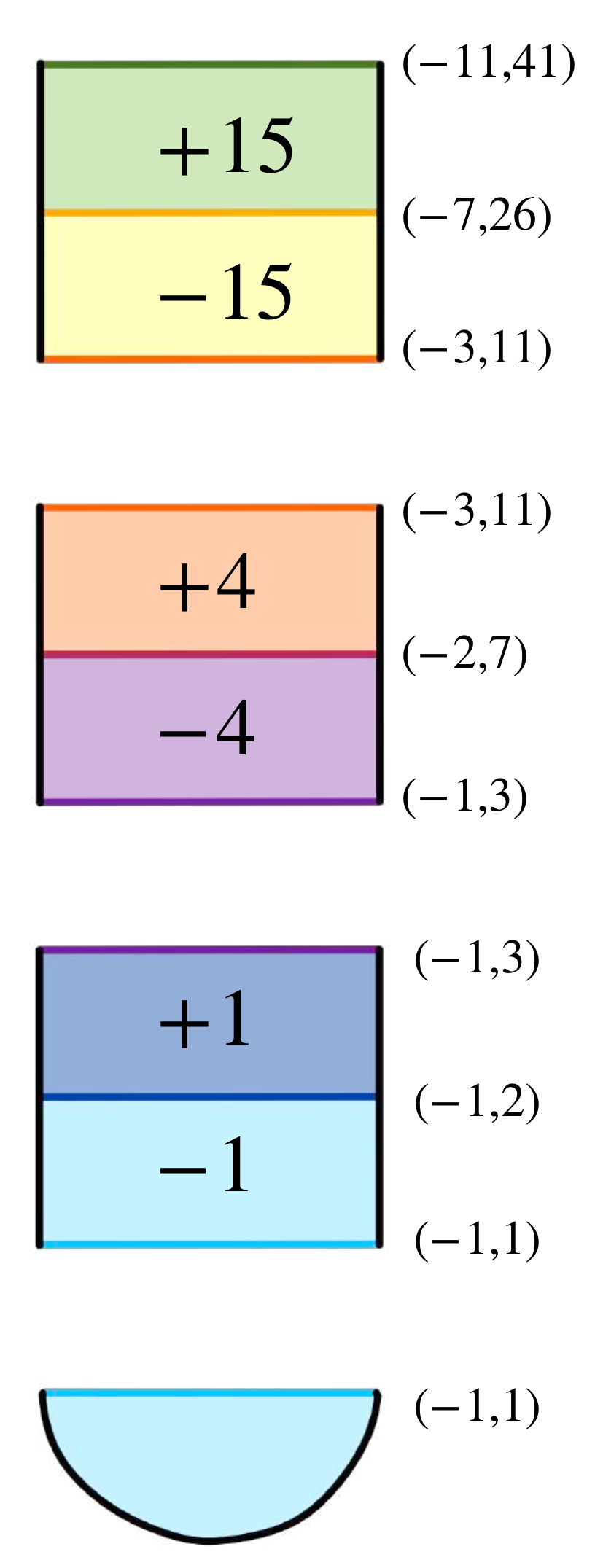}
  \caption{$L(56,15)$}
  \label{blocksL5615}
\end{subfigure}
\begin{subfigure}[t]{.5\textwidth}
  \centering
 \includegraphics[scale=0.5]{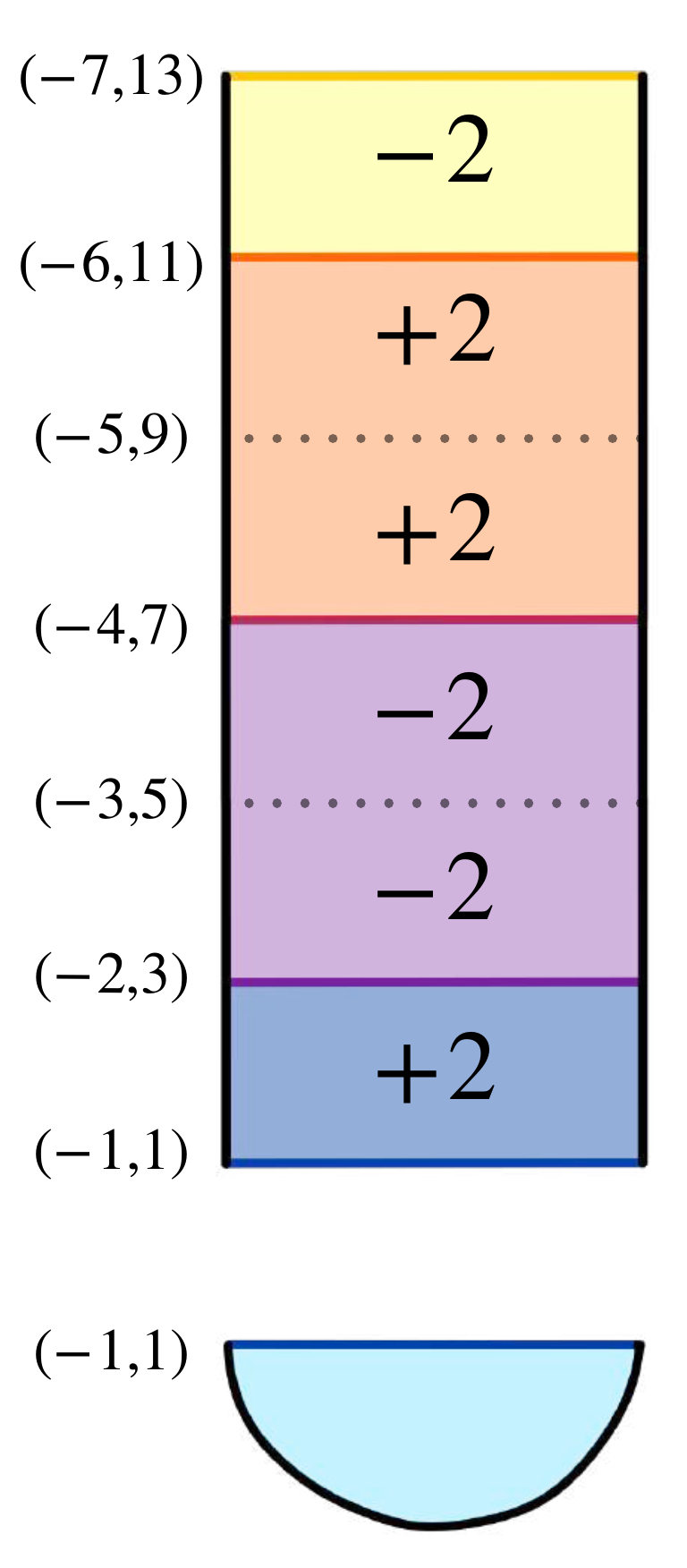}
  \caption{$L(28,15)$}
  \label{blocksL2815}
\end{subfigure}
\caption{Compatible choice of signs for a covering map.}
\label{L5615L2815}
\end{figure}

\begin{figure}[ht!]
\centering
\begin{subfigure}[t]{.45\textwidth}
  \centering
 \includegraphics[scale=0.45]{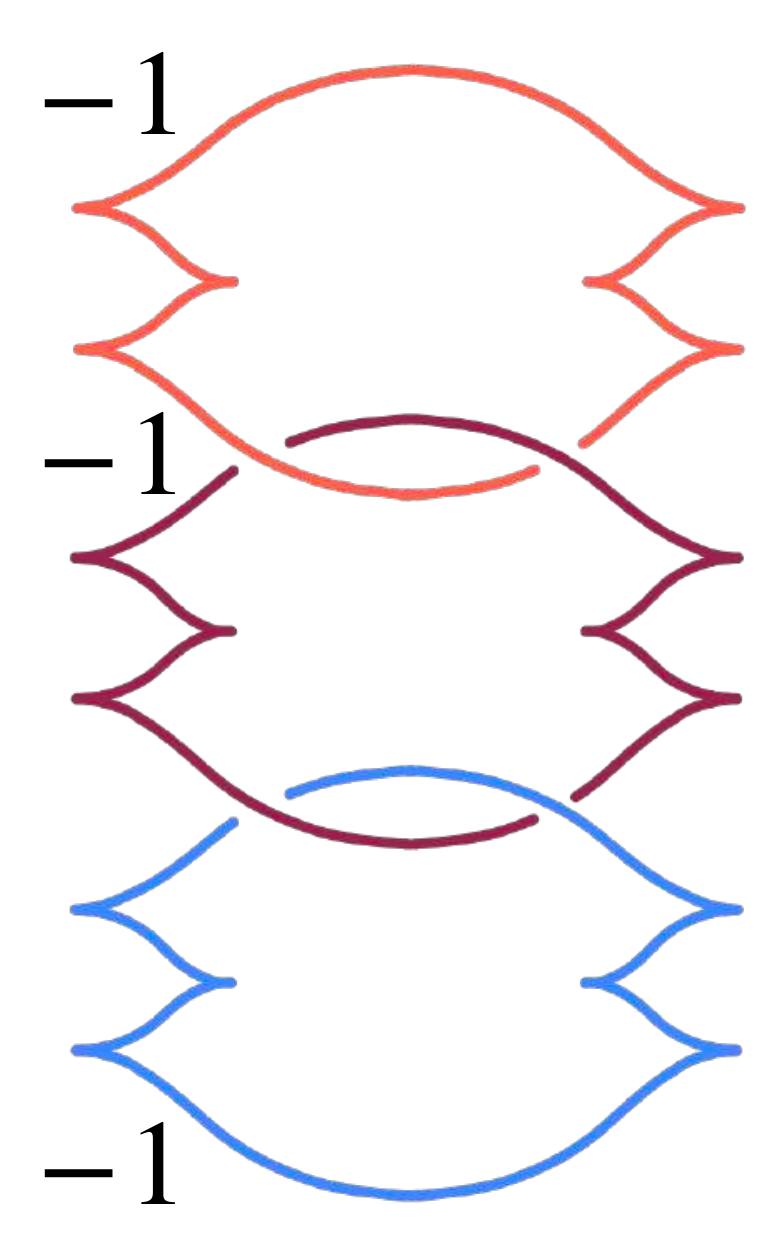}
\caption{$L(56,15)$}
  \label{L5615}
\end{subfigure}
\begin{subfigure}[t]{.45\textwidth}
  \centering
 \includegraphics[scale=0.5]{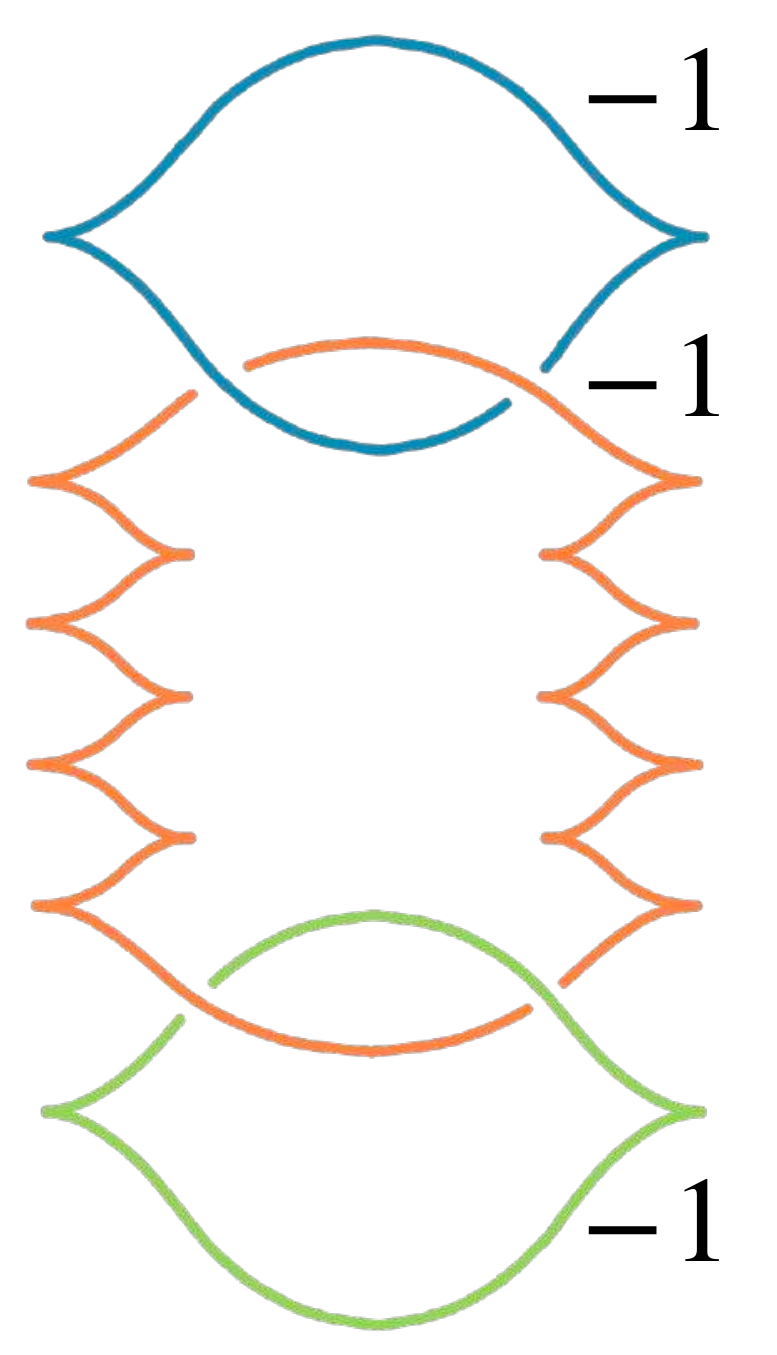}
  \caption{$L(28,15)$}
  \label{L2815}
\end{subfigure}
\caption{Contact surgery producing lens spaces.}
\label{L5615L2815links}
\end{figure}

A compatible case is illustrated for example by Figure \ref{L5615L2815},
which represents the double cover
\[L(28,15)\to L(56,15),\]
where the contact structures on the two lens spaces are specified by Figure \ref{L5615L2815links}.

If we look at the lattice of coverings of $L(56,15)$ we see that this contact structure $\xi$ (Figure \ref{L5615}) lifts to an overtwisted one along some (at least one) covering maps:

\[
\xymatrix{  & (L(7,1),{\color{blue}\xi_{?}})\ar[r]^-{2:1} & (L(14,1),{\color{blue}\xi_{?}}) \ar[r]^-{2:1}& (L(28,15),{\color{blue}\xi_{vot}})\ar[rd]^-{2:1} \\
(S^3,{\color{red}\xi_{ot}}) \ar[ru]^-{7:1}\ar[rd]_-{2:1}& & & &  (L(56,15),{\color{blue}\xi_{vot}}) \\
 & (L(2,1),{\color{red}\xi_{ot}})\ar[r]_-{2:1})\ar[ruu]_-{7:1} & (L(4,3),{\color{red}\xi_{ot}}) \ar[r]_-{2:1}\ar[ruu]_-{7:1}& (L(8,7),{\color{red}\xi_{ot}})\ar[ru]_-{7:1} 
 }\]
 
Therefore we cannot apply directly the criterion of previous section to conclude that the Stein fillings of $(L(56,15),\xi)$ are simply connected. By the fact that lifting $\xi$ to $L(8,7)$ results in an overtwisted structure, we get that the kernel of $i_*$ cannot be contained in $\Z/8\Z$, where
\[i:L(56,15)\hookrightarrow X\]
is the inclusion of the boundary of any Stein filling $X$. We have that
\[\pi_1(X)=\frac{\Z/56\Z}{\ker i_*},\]
so the possibilities are:
\begin{itemize}
\item $\ker i_*=\Z/7\Z$, which gives $\pi_1(X)=\Z/8\Z$,
\item $\ker i_*=\Z/14\Z$, which gives $\pi_1(X)=\Z/4\Z$,
\item $\ker i_*=\Z/28\Z$, which gives $\pi_1(X)=\Z/2\Z$,
\item $\ker i_*=\Z/56\Z$, which gives $\pi_1(X)=1$. The following proposition proves that  
\end{itemize}

\begin{prop} \label{pi=1}
$\pi_1(X)=1$.
\end{prop}

\begin{prf} 
Consider the Stein filling $X_{\Lambda}$ of $(L(56,15),\xi)$ described by the diagram of Figure \ref{L5615}. We want to compute the $d_3$ invariant of the contact structure on the boundary:
\[d_3(\xi)=\frac{1}{4}(c_1(X_{\Lambda})^2 -3\sigma(X_{\Lambda})-2\chi(X_{\Lambda})).\]
The first Chern class $c_1(X_{\Lambda})$ is zero, because it evaluates as $\rot_i=0$ on the three generators of $H_2(X_{\Lambda})$. Moreover, $\sigma(X_{\Lambda})=-3$ and $\chi(X_{\Lambda})=4$. Therefore
\[d_3(\xi)=\frac{1}{4}.\]
Notice that $c_1(\xi)=0$ because it is the restriction of $c_1(X_{\Lambda})$, which is 0 itself. Being any contact structure on a lens space planar (\cite{schonenberger}), we can apply \cite[Corollary 1.5]{c1planar} and conclude that any Stein filling of $(L(56,15),\xi)$ has vanishing $c_1$. 

We want to compute $d_3(\xi)$ using the Stein filling $X$. For what we have just said $c_1(X)^2=0$ and we also have $\sigma(X)=1-\chi(X)$. So:
\[\frac{1}{4}=d_3(\xi)=\frac{1}{4}(c_1(X)^2 -3\sigma(X)-2\chi(X))=\frac{1}{4}(-3+\chi(X)).\]
This tells us that 
\[\fbox{$\chi(X)=4.$}\]

Now we analyze the possibilities for its fundamental group case by case.
\begin{itemize} 
\item[i)] Suppose that $\pi_1(X)=\Z/8\Z.$ Then we pass to the universal covering $\widetilde{X}\to X$, of degree 8, whose boundary is the (connected) degree-8 covering $L(7,1)\to L(56,15)$. By Theorem \ref{pi1b2} we have
\[\chi(X)\leq \frac{1+\length(7/1)}{8}=\frac{1+1}{8}=\frac{1}{4},\]
which is impossible. So $\pi_1(X)\neq \Z/8\Z$.

\item[ii)] If $\pi_1(X)=\Z/4\Z$, we pass to the universal covering and since $\length(14/1)=1$, we get $\chi(X)\leq 1/2$. This is not possible, hence $\pi_1(X)\neq \Z/4\Z$. 

\item[iii)] Again, we take the universal covering $\widetilde{X}\to X$, of degree 2, whose boundary is the (connected) degree-2 covering $L(28,15)\to L(56,15)$. By Theorem \ref{pi1b2}, we have
\[\chi(X)\leq \frac{1+\length(28/15)}{2}=\frac{1+3}{2}=2\]
and hence $\chi(X)\leq 2$, which is not possible. Hence $\pi_1(X)\neq \Z/2\Z$.

\item[iv)] We conclude that any Stein filling of $(L(56,15),\xi)$ is simply connected. (Note that by \cite[Theorem 1.3]{menke} we already know that in fact there is a unique filling obtained by attaching three 2-handles to $B^4$ along the link of Figure \ref{L5615}). 
\end{itemize}
\end{prf}
 
The result proved in Theorem \ref{pi1b2} is that somehow for a Stein filling $X$ of $(L(p,q),\xi_{vo})$ "the bigger $\pi_1(X)$ is, the smaller its Euler characteristic is forced to be". Of course this is in general spoiled by the quantity $l'$, appearing in the statement, which depends on the numbers $p/d$ and $q$ (one should first reduce $q$ modulo $p/d$, in case it were bigger).

On the other hand, if $p$ is small, then by Theorem \ref{thmcovering} we have a bigger chance of finding coverings of $L(p,q)$ which are overtwisted, and hence apply our criterion to bound the cardinality of $\pi_1(X)$. For this reason, in the first version of this paper we expected that every Stein filling of $(L(p,q),\xi_{vo})$ is simply connected.

\noindent \textbf{However}, Marco Golla found an easy way to produce Stein fillings of virtually overtwisted structures on lens spaces that are not simply connected: for example, let 
\[-\frac{p}{q}=[-4,-2n,-4], \;n>1, \]
and consider the Legendrian representative of the 3-components link associated to this continued fraction expansion where the first and third components have rotation number +2, while the middle one has rotation number 0. The fillings of this virtually overtwisted structure are understood thanks to the work of \cite{menke} and \cite{mcduff}. In particular, there is a filling which is obtained from a boundary connected sum of two rational homology balls with $\pi_1=\Z/2\Z$ (corresponding to the two $-4$) by attaching a single Weinstein 2-handle (corresponding to the central $-2n$): this handle attachment does not kill the whole $\Z/2\Z\ast \Z/2\Z$, resulting in a non simply-connected filling.

\bibliographystyle{alpha} 

\bibliography{Topological_constraints.bib}
\thispagestyle{plain}

\Addresses

\end{document}